\documentclass[11pt, one side,emlines]{amsart}
\usepackage{amssymb,latexsym,xy,eucal,mathrsfs, graphicx, tikz, amsmath, caption, tabularx}
\usepackage{ circuitikz }
\usepackage{algorithm}
\usepackage{algorithmic}
\textwidth=15.2cm \textheight=24cm \theoremstyle{plain}
\newtheorem{lem}{Lemma}[section]
\newtheorem{thm}[lem]{Theorem}

\newtheorem{cor}[lem]{Corollary}

\theoremstyle{definition}
\newtheorem{definition}[lem]{Definition}

\newtheorem{prep}{Preposition}
\newtheorem{note}{Note}

\numberwithin{equation}{section}  \voffset
-55truept \hoffset -20truept
\begin{document}
	\baselineskip 15truept
	\title[Construction and Conditions for CISTs in Hypercubes and Bipartite Graphs]{Construction and Conditions for Completely Independent Spanning Trees in Hypercubes and Regular Bipartite Graphs}
	\maketitle 
		\begin{center} 
		
			R. Barabde,	S. A. Mane, S. A. Kandekar$^*$\\  
				{\small Department of Mathematics,\\ Savitribai Phule Pune University, Pune-411007, India.}\\
				\email{\emph{rushikeshbarbde@gmail.com; manesmruti@yahoo.com; smitakandekar54@gmail.com }} 
			\end{center}
	\begin{abstract}

	  		A set of \( k \) spanning trees in a graph \( G \) is called a set of \textit{completely independent spanning trees (CISTs)} if, for every pair of vertices \( x \) and \( y \), the paths connecting \( x \) and \( y \) across different trees do not share any vertices or edges, except for \( x \) and \( y \) themselves. Hasunuma conjectured that every \(2k\)-connected graph contains exactly \(k\) completely independent spanning trees (CISTs). However, Pétérfalvi disproved this conjecture.
	  	
	  	When \( k = 2 \), the two CISTs are called a \textit{dual-CIST}. It has been shown that determining whether a graph can have \( k \) CISTs is an NP-complete problem, even when \( k = 2 \). In $2017$, Darties et al. raised the question of whether the $6-$dimensional hypercube \( Q_6 \) can have three completely independent spanning trees (CISTs). This paper provides an answer to that question.
	  	
	  	In this paper, we first present a necessary condition for \( k \)-regular, \( k \)-connected bipartite graphs to have \( \left\lfloor \frac{k}{2} \right\rfloor \) CISTs. 	We also investigate that the hypercube of dimension \( n \) cannot have \( \frac{n}{2} \) CISTs, which means Hasunuma's conjecture does not hold for the hypercube \( Q_n \) when \( n \) is an even integer \(2 < n \leq 10^7 \), except when \(n = 2^r\) and \( n \in \{161038, 215326, 2568226, 3020626, 7866046, 9115426 \} \). This result also resolves a question posed by Darties et al. 
	  	
	  	 The construction of multiple CISTs on the underlying graph of a network has practical applications in ensuring the fault tolerance of data transmission. In this context, we also provide a construction for three completely independent spanning trees in the hypercube \(Q_n\) for \(n \geq 7\). Our results show that Hasunuma's conjecture holds for odd integer \(n = 7\) in \(Q_n\), but does not hold for even integer \(n = 6\).
	  	
	  	\end{abstract}
	
	\section{Introduction} 
The network topology is represented as a graph, where vertices correspond to nodes and edges symbolize direct physical connections between these nodes. Interconnection networks are commonly represented as undirected simple graphs \( G = (V, E) \). 
	
	A tree \( T \) is called a spanning tree of a graph \( G \) if \( V(T) = V(G) \). A vertex \( v \in V(T) \) is termed internal if its degree in \( T \), denoted \( d_T(v) \), is greater than or equal to 2. Two spanning trees \( T_1 \) and \( T_2 \) in \( G \) are said to be edge-disjoint if their edge sets are disjoint, i.e., \( E(T_1) \cap E(T_2) = \emptyset \). For a given tree \( T \) and a pair of vertices \( u, v \in V(T) \), let \( P_T(u,v) \) represent the set of vertices along the unique path between \( u \) and \( v \) in \( T \). Two spanning trees \( T_1 \) and \( T_2 \) are internally vertex-disjoint if for any pair of vertices \( u \) and \( v \) in \( V(G) \), the intersection of the paths between \( u \) and \( v \) in \( T_1 \) and \( T_2 \) is restricted to the endpoints, i.e., \( P_{T_1}(u,v) \cap P_{T_2}(u,v) = \{u, v\} \). Finally, a collection of spanning trees \( T_1, T_2, \ldots, T_k \) of \( G \) is termed \textit{pairwise edge-disjoint} if the edge sets of any two trees are disjoint, and \textit{pairwise internally vertex-disjoint} if the paths between any two vertices in different trees intersect only at the endpoints. Such a collection of spanning trees is called \textit{completely independent spanning trees} (CISTs for short).
	
	 The concept of CISTs is inspired by the early work of Hasunuma \cite{h1}. Hasunuma conjectured that any \( 2k \)-connected graph contains \( k \) CISTs. However, this conjecture was disproven by Pétrfalvi \cite{pt}. The notion of CISTs emerged from the study of routing, reliability, and fault tolerance in parallel and distributed systems. Various authors have since provided necessary and sufficient conditions for the existence of CISTs \cite{C1,C2,C3,ta}. Further work on CISTs can be found in \cite{ c3, h1, sm, p3, p6}. Hasunuma \cite{h2} observed that determining the existence of a dual-CIST is an NP-complete problem
	
	Pai et al. \cite{p1} gave a necessary condition for \( k \)-regular, \( k \)-connected graphs to admit \( \left\lfloor \frac{k}{2} \right\rfloor \) CISTs. 
	
	\begin{thm}\cite{p1} Let \( G \) be a \( k \)-connected, \( k \)-regular graph. If \( G \) admits \( \lfloor \frac{k}{2} \rfloor \) CISTs, then the following condition holds
		
		\[
		\left \lceil \frac{{|V(G)|} - 2}{\lceil \frac{k}{2} \rceil} \right\rceil \leq \left\lfloor \frac{{|V(G)|}}{\lfloor \frac{k}{2} \rfloor} \right\rfloor
		\]
	\end{thm}

	They also demonstrated that this conjecture does not hold for hypercubes of dimensions \( n \in \{10, 12, 14, 20, 22, 24, 26, 28, 30\} \). In this paper, we provide a necessary condition for \( k \)-regular, \( k \)-connected bipartite graphs to have \( \left\lfloor \frac{k}{2} \right\rfloor \) CISTs. 
	We also investigate that Hasunuma's conjecture does not hold for hypercube $Q_n$ with all even integers  \(n, 2 < n \leq 10^7 \), except when \(n = 2^r\) and \( n \in \{161038, 215326, 2568226, 3020626, 7866046, 9115426 \} \). 
This result also resolves a question posed by Darties et al. \cite{dar} in $2017$, who asked whether the $6-$dimensional hypercube \( Q_6 \) could have three completely independent spanning trees.
	
	The construction of multiple CISTs on the underlying graph of a network has practical applications in ensuring the fault tolerance of data transmission. Various methods for constructing three CISTs in non-bipartite hypercube variants of dimension $6$ have been established for several hypercube-related graphs, including crossed cubes, Möbius cubes, and shuffle cubes. \cite{c3, p3, p6}. In this paper, we present a construction for three completely independent spanning trees in \( Q_n \) for \( n \geq 7 \), thereby proving that Hasunuma's conjecture holds for odd integer \(n = 7\) in \(Q_n\), but does not hold for even integer \(n = 6\).
	
	For additional terminology and notation, refer to \cite{we}.

		\section{\textbf{Preliminaries}}

		Interconnection networks have been widely studied recently. The architecture of an interconnection networks are usually denoted as a graph $G$. A graph $G$ is an ordered pair $(V(G), E(G) )$ comprising a set of vertices or nodes $V(G)$ together with a set $E(G)$ of edges which are two element subsets of $V(G).$ Many useful topologies have been proposed to balance performance and cost parameters. Among them, the hypercube denoted by $Q_n$ is one of the most popular topologies and has been studied for parallel networks. 
		
	The \( n \)-dimensional hypercube, \( Q_n \), is a graph with \( 2^n \) vertices, each labeled by an \( n \)-bit binary string. Specifically, \( Q_1 \) is the graph \( K_2 \) with vertex set \( \{0, 1\} \). For \( n \geq 2 \), \( Q_n \) can be recursively constructed from two copies of \( Q_{n-1} \), denoted \( Q^0_{n-1} \) and \( Q^1_{n-1} \), by adding \( 2^{n-1} \) edges between them. The vertex set \( V(Q^0_{n-1}) \) consists of strings starting with \( 0 \), and \( V(Q^1_{n-1}) \) consists of strings starting with \( 1 \). A vertex \( u = 0u_2 \ldots u_n \) in \( Q^0_{n-1} \) is adjacent to a vertex \( v = 1v_2 \ldots v_n \) in \( Q^1_{n-1} \) if \( u_i = v_i \) for all \( i \geq 2 \); the edge \( \langle u, v \rangle \) is called a hypercube edge, and we write \( v = u^h \).
	
	Furthermore, the vertex set of \( Q_n \) can be divided into two sets \( X \) and \( Y \) where, \( X \) consists of vertices whose labels contain an even number of 1's 
	and \( Y \) consists of vertices whose labels contain an odd number of 1's.
	
	By the definition of the hypercube, both \( X \) and \( Y \) are independent sets, and they partition the vertex set of \( Q_n \). Hence, \( Q_n \) is a bipartite graph with partition sets \( X \) and \( Y \).
	
	The edge set of \( Q_n \) is denoted by \( E_n = \{\langle u, u^h \rangle : u \in V(Q^0_{n-1})\} \). See Figure $1$.\\
\vspace*{-0.5 cm}		
\begin{figure}[!ht]
	\centering
	\resizebox{0.7\textwidth}{!}{%
		\begin{circuitikz}
			\tikzstyle{every node}=[font=\Large]
			\node at (1.25,5.75) [circ] {};
			\node at (1.25,8.25) [circ] {};
			\node at (2.5,8.25) [circ] {};
			\node at (2.5,5.75) [circ] {};
			\node at (5,8.25) [circ] {};
			\node at (5,5.75) [circ] {};
			\node at (6.25,8.25) [circ] {};
			\node at (6.25,5.75) [circ] {};
			\node at (8.75,8.25) [circ] {};
			\node at (8.75,5.75) [circ] {};
			\node at (10,8.25) [circ] {};
			\node at (12.5,8.25) [circ] {};
			\node at (10,5.75) [circ] {};
			\node at (12.5,5.75) [circ] {};
			\draw [ line width=1pt](1.25,8.25) to[short] (1.25,5.75);
			\draw [ line width=1pt](2.5,8.25) to[short] (2.5,5.75);
			\draw [ line width=1pt](2.5,8.25) to[short] (5,8.25);
			\draw [ line width=1pt](5,8.25) to[short] (5,5.75);
			\draw [ line width=1pt](2.5,5.75) to[short] (5,5.75);
			\draw [ line width=1pt](6.25,8.25) to[short] (6.25,5.75);
			\draw [ line width=1pt](6.25,8.25) to[short] (8.75,8.25);
			\draw [ line width=1pt](8.75,8.25) to[short] (8.75,5.75);
			\draw [ line width=1pt](6.25,5.75) to[short] (8.75,5.75);
			\draw [ line width=1pt](10,8.25) to[short] (10,5.75);
			\draw [ line width=1pt](10,8.25) to[short] (12.5,8.25);
			\draw [ line width=1pt](12.5,8.25) to[short] (12.5,5.75);
			\draw [ line width=1pt](10,5.75) to[short] (12.5,5.75);
			\draw [line width=1pt, short] (6.25,8.25) -- (8.25,8.25);
			\draw [line width=1pt, short] (6.25,8.25) .. controls (8.25,9) and (8.5,8.75) .. (10,8.25);
			\draw [line width=1pt, short] (8.75,8.25) .. controls (10.5,9) and (11,9) .. (12.5,8.25);
			\draw [line width=1pt, short] (6.25,5.75) .. controls (8,5) and (8.5,5) .. (10,5.75);
			\draw [line width=1pt, short] (8.75,5.75) .. controls (10.5,5) and (11.25,5) .. (12.5,5.75);
			\node [font=\Large] at (1,4.7) {Q};
			\node [font=\normalsize] at (1.25,4.45) {1};
			\node [font=\Large] at (3.75,4.7) {Q};
			\node [font=\Large] at (9.25,4.7) {Q};
			\node [font=\normalsize] at (4,4.45) {2};
			\node [font=\normalsize] at (9.5,4.45) {3};
			\node [font=\normalsize] at (1,5.5) {0};
			\node [font=\normalsize] at (1,8.25) {1};
			\node [font=\normalsize] at (2.25,5.5) {00};
			\node [font=\normalsize] at (2.25,8.5) {01};
			\node [font=\normalsize] at (4.75,5.5) {10};
			\node [font=\normalsize] at (5,8.5) {11};
			\node [font=\normalsize] at (6,5.5) {000};
			\node [font=\normalsize] at (8.25,5.5) {010};
			\node [font=\normalsize] at (6,8.5) {001};
			\node [font=\normalsize] at (8.25,8.5) {011};
			\node [font=\normalsize] at (10.25,5.5) {100};
			\node [font=\normalsize] at (12.75,5.5) {110};
			\node [font=\normalsize] at (12.75,8.5) {111};
			\node [font=\normalsize] at (10.25,8.5) {101};
		\end{circuitikz}
	}%
	\vspace*{-0.3 cm}	
	\caption{Hypercubes of dimension 1, 2, and 3}
\end{figure}
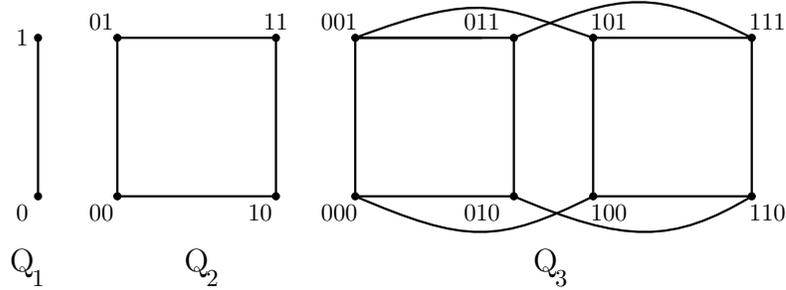

\section {\textbf{Necessary condition for \( k \)-regular, \( k \)-connected bipartite graphs to have \( \left\lfloor \frac{k}{2} \right\rfloor \) CISTs.} }

In this section, we will provide a necessary condition for \( k \)-regular, \( k \)-connected bipartite graphs to have \( \left\lfloor \frac{k}{2} \right\rfloor \) CISTs. Note that if \( k > 0 \), then a \( k \)-regular bipartite graph has the same number of vertices in each partite set \cite{we}. 

	\begin{thm} Let \( G \) be a \( k \)-regular, \( k \)-connected bipartite graph. If \( G \) contains \( \lfloor \frac{k}{2} \rfloor \) CISTs, then the following inequality holds
		
		\[
		\left \lceil \frac{\frac{|V(G)|}{2} - 1}{\lceil \frac{k}{2} \rceil} \right\rceil \leq \left\lfloor \frac{\frac{|V(G)|}{2}}{\lfloor \frac{k}{2} \rfloor} \right\rfloor
		\]
\end{thm}

\begin{proof}

Let \( G \) be a \( k \)-regular, \( k \)-connected bipartite graph with a vertex partition \( V(G) = X \cup Y \), where \( X \) and \( Y \) are disjoint independent sets, and \( |X| = |Y| = \frac{|V(G)|}{2} = m \).

Suppose \( G \) contains \( \lfloor \frac{k}{2} \rfloor \) CISTs. Let T be a spanning tree of G among these CISTs such that the number of internal vertices of T in set X is minimum. Let S be the set of internal vertices of T. Then, \( S = S_1 \cup S_2 \), where \( S_1 \subseteq X \) and \( S_2 \subseteq Y \), and both \( S_1 \) and \( S_2 \) are non-empty.

Clearly, all the vertices in \( X \setminus S_1 \) are leaves in \( T \). Therefore, when we remove all the vertices in $X\backslash S1$ then the induced subgraph of T on the vertices $S1 \cup Y$ is connected. This subgraph has \( |Y \cup S_1| - 1 \) edges, and every edge has one endpoint in \( S_1 \). This implies that the total degree  \( |Y \cup S_1| - 1 \) of the vertices in \( S_1 \) is used in forming the tree \( T \).

By the definition of a CIST, every vertex of $G$ can be a non-leaf in at most one CIST and a leaf in others. Therefore, for the \( \lfloor \frac{k}{2} \rfloor \) CISTs in the \( k \)-regular graph \( G \), the degree of any non-leaf vertex in a CIST is at most \( k - \lfloor \frac{k}{2} \rfloor + 1 \). This implies:

\[
|Y \cup S_1| - 1 \leq |S_1| \cdot \left( k - \left\lfloor \frac{k}{2} \right\rfloor + 1 \right)
\]

Since Y and S1 are disjoint sets,

\[
|Y| + |S_1| - 1 \leq |S_1| \cdot \left( \left\lceil \frac{k}{2} \right\rceil \right) + |S_1|
\]

Since \( |Y| = m \), this simplifies to

\[
m - 1 \leq |S_1| \cdot \left\lceil \frac{k}{2} \right\rceil
\]

Thus,

\[
|S_1| \geq \frac{m - 1}{\left\lceil \frac{k}{2} \right\rceil}
\]

Since \( |S_1| \) is an integer, we can therefore write

\[
|S_1| \geq \left\lceil \frac{m - 1}{\left\lceil \frac{k}{2} \right\rceil} \right\rceil \tag{1}
\]

Next, as every vertex of $G$ can serve as a non-leaf in at most one CIST and \( S_1 \) is the minimal subset of \( X \) containing non-leaf vertices, \( S_1 \) contains at most \( \frac{|X|}{\lfloor \frac{k}{2} \rfloor} \) vertices. This gives,
\[
|S_1| \leq \left\lfloor \frac{m}{\lfloor \frac{k}{2} \rfloor} \right\rfloor \tag{2}
\]

Combining inequalities (1) and (2), we obtain:
\[
\left\lceil \frac{m - 1}{\left\lceil \frac{k}{2} \right\rceil} \right\rceil \leq \left\lfloor \frac{m}{\lfloor \frac{k}{2} \rfloor} \right\rfloor
\]

Thus, we have the desired inequality:
\[
\left \lceil \frac{\frac{|V(G)|}{2} - 1}{\lceil \frac{k}{2} \rceil} \right\rceil \leq \left\lfloor \frac{\frac{|V(G)|}{2}}{\lfloor \frac{k}{2} \rfloor} \right\rfloor \tag{3}
\]

\end{proof}

\begin{note}
	For the 6-dimensional hypercube \( Q_6 \), if we substitute the required values into the left-hand side (LHS) and right-hand side (RHS) of the equation given by Theorem 1.1, we get LHS = RHS = 21. However, if we substitute these values into the equation number $3$ given by Theorem 3.1, we get LHS = 11 and RHS = 10. Therefore, by Theorem 3.1, it is proven that \( Q_6 \) does not admit three completely independent spanning trees (CISTs). This answers the question posed by Darties et al\cite{dar}.
	\end{note}

\section {\textbf{Hasunuma's conjecture in the hypercube of even dimension.} }

In this section, we will prove that there do not exist \( \frac{n}{2} \) completely independent spanning trees in the hypercube \( Q_n \) when \( n \) is an even integer such that \( 2 < n \leq 10^7 \), except when \( n = 2^r \)  and \( n \in \{161038, 215326, 2568226, 3020626, 7866046, 9115426 \} \). Before proceeding with the proof, we first state and prove a lemma, followed by an algorithm, and then we present the main proof.

\begin{lem}
	
	Let \( x \) and \( k \) be positive integers such that \( k \nmid x \) and \( k \nmid x-1 \). Then,
	
	\[
	\left \lceil \frac{x-1}{k} \right\rceil > \left\lfloor \frac{x}{k} \right\rfloor.
	\]
	
\end{lem}

\begin{proof} Suppose, for the sake of contradiction, that
	
	\[
	\left \lceil \frac{x-1}{k} \right\rceil \leq \left\lfloor \frac{x}{k} \right\rfloor.
	\]
	
	Since \( k \nmid x \) and \( k \nmid x-1 \), both \( \frac{x-1}{k} \) and \( \frac{x}{k} \) are non-integer values.
	
	Let \( \left\lfloor \frac{x}{k} \right\rfloor = a \), which implies \( \frac{x}{k} > a \). Also, \( \left \lceil \frac{x-1}{k} \right\rceil \leq a \), which implies \( \frac{x-1}{k} < a \). 
	
	By combining these inequalities, we have
	
	\[
	\frac{x-1}{k} < a < \frac{x}{k},
	\]
	
	which implies
	\[
	x - 1 < ak < x.
	\]
	
	This leads to a contradiction because \( x \) and \( x-1 \) are consecutive integers, while \( ak \) is an integer. Therefore, the assumption that \( \left \lceil \frac{x-1}{k} \right\rceil \leq \left\lfloor \frac{x}{k} \right\rfloor \) must be false.
	
	Thus, we conclude that
	
	\[
	\left \lceil \frac{x-1}{k} \right\rceil > \left\lfloor \frac{x}{k} \right\rfloor \tag{4}.
	\]
	
\end{proof}
Before we proceed with the main proof, we present the following algorithm that identifies all even integers \( m \leq \text{limit} \), where \( \text{limit} \) is a positive integer, such that the expression

\[
\frac{2^{m-1} - 1}{\frac{m}{2}}
\]

is an integer. For convenience, we take \( m \leq 10^7 \).

\begin{algorithm}
	\noindent
\caption{Finding Valid Even Values of \( m \)}\label{alg:find_valid_m}
\begin{flushleft}
	\textbf{Input:} limit (\textit{int}) -- Upper limit for \( m \)\\
\textbf{Output:} List of valid even \( m \)\\
valid\_m \(\gets\) []\\
{\( m \gets 2, 4, \ldots, \text{limit} \)} {Iterate over even values of \( m \)}\\
numerator \(\gets 2^{m-1} - 1\)\\
denominator \(\gets m / 2\)\\
If, {numerator \(\bmod\) denominator \( = 0\)}\\
Append \( m \) to valid\_m
	
\end{flushleft}
	
\end{algorithm}

From Algorithm \ref{alg:find_valid_m}, the even integers \( m \in \{161038, 215326, 2568226, 3020626, 7866046, 9115426 \} \) such that \( 2 < m \leq 10^7 \), are the only even positive integers for which

\[
\frac{2^{m-1}-1}{\frac{m}{2}}
\]

is an integer.

\newpage
\begin{prep}

 Let \( n \) and \( r \) be integers such that \( n \) is an even integer with \( 2 < n \leq 10^7 \). Then, there do not exist \( \frac{n}{2} \) completely independent spanning trees in the hypercube \( Q_n \), except when \( n = 2^r \) and \( n \in \{161038, 215326, 2568226, 3020626, 7866046, 9115426 \} \).
\end{prep}

\begin{proof}
	
	By Lemma 4.1, we have
	
\[
\left \lceil \frac{x-1}{k} \right\rceil > \left\lfloor \frac{x}{k} \right\rfloor,
\]
for positive integers x and k such that \( k \nmid x \) and \( k \nmid x-1. \)
we substitute \( x = \frac{|V(Q_n)|}{2} = \frac{2^n}{2} = 2^{n-1} \) for all even \( n = 2k > 2 \). Then, for \( n \neq 2^r \), the factors of \( k \) will always contain at least one odd integer, and hence the right-hand side (RHS) of this inequality is always a non-integer value. 

If we substitute these values of $x$ and $k$ into the left-hand side (LHS) of this inequality, then according to the algorithm above, the LHS is a non-integer  except for the values \( n \in \{161038, 215326, 2568226, 3020626, 7866046, 9115426 \}. \) Thus, since both the LHS and RHS are non-integers, the inequality holds true for above value of x and k, i.e.

\[
\left \lceil \frac{\frac{|V(Q_n)|}{2} - 1}{\lceil\frac{n}{2}\rceil}\right\rceil > \left\lfloor \frac{\frac{|V(Q_n)|}{2}}{\lfloor\frac{n}{2}\rfloor}\right\rfloor.
\] 

Thus, by Theorem 3.1, for an even integer \( n \) with \( 2 < n \leq 10^7, Q_n \) does not admit \(  \frac{n}{2} \) CISTs, except when \( n = 2^r \) and \( n \in \{161038, 215326, 2568226, 3020626, 7866046, 9115426 \} \). 

\end{proof}

\begin{note}
 In the algorithm, we set the limit to \( m = 10^7 \), but we can increase the limit to obtain more values of the dimension \( n \) of hypercubes \( Q_n \) as required. Also, for \( n=2^r\) with \( 2 < n \leq 10^7\) and \(n \in \{161038, 215326, 2568226, 3020626, 7866046, 9115426\} \), the inequality (1) holds true. However, this does not guarantee that such \( Q_n \) will always admit \( \left\lfloor \frac{n}{2} \right\rfloor \) completely independent spanning trees (CISTs).
\end{note}

 \vspace*{0.5 cm}	
		\section {\textbf{Existence of three CISTs in the hypercube \( Q_7 \)} }
		After showing that Hasunuma's conjecture is disproved for even dimensions \( n \) of \( Q_n \), except for some specific values of \( n \) as discussed above. 
		In this section, we will now demonstrate that there exist \( \frac{n}{2} \) completely independent spanning trees (CISTs) in \( Q_n \) for \( n = 7 \). Hence, Hasunuma's conjecture holds true for the odd integer \( n = 7 \).
		  Here, among these spanning trees, two of these trees have diameters of length  \(15\) and \(18\), respectively, and the third tree has a diameter of length \(17\).		
		
		
		For clarity, we represent the vertices of \( Q_7 \) using decimal notations.
		
		Let \( InV(T) \) denote the set of internal vertices of the tree \( T \). A vertex \( v \in V(T) \) is considered internal if its degree \( d_{T}(v) \) is at least 2.
		
		We first enumerate three Complete Independent Spanning Trees (CISTs) in \( Q_7 \). Refer to figures 2, 3, and 4.		
		
\vspace*{-0.7 cm}
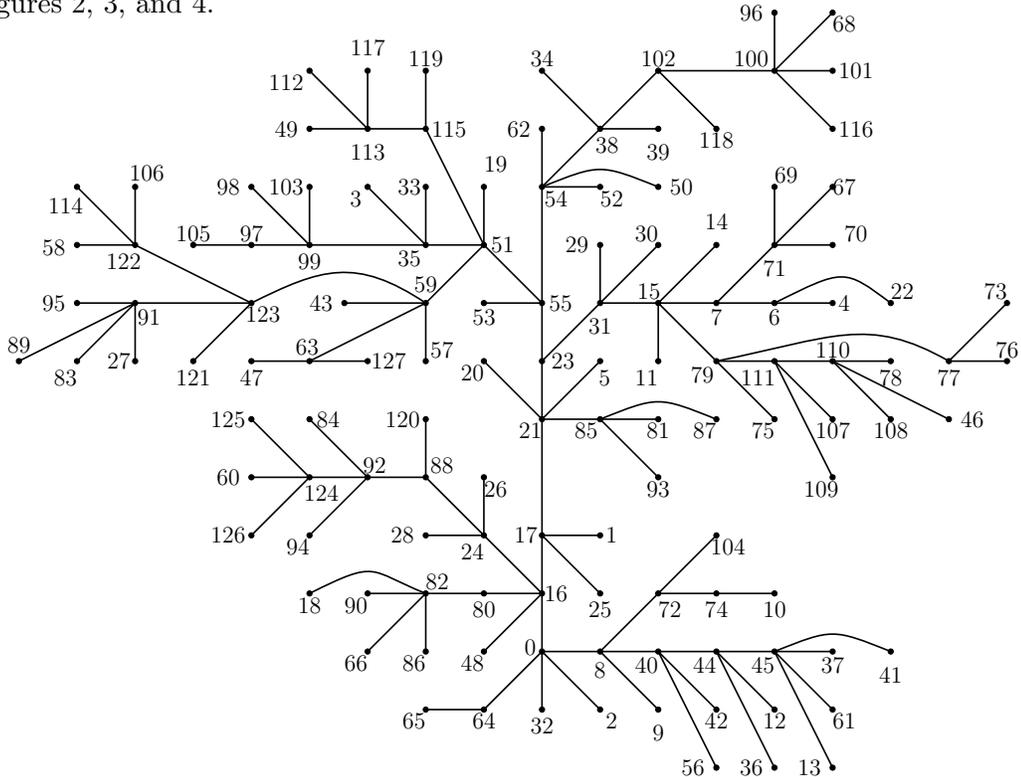
\begin{figure}[!ht]
	\centering
	\resizebox{0.9\textwidth}{!}{%
		\begin{circuitikz}
			\tikzstyle{every node}=[font=\Large]
			\node at (6.25,8.75) [circ] {};
			\node at (5,8.75) [circ] {};
			\node at (3.75,8.75) [circ] {};
			\node at (7.5,8.75) [circ] {};
			\node at (2.5,10) [circ] {};
			\node at (3.75,10) [circ] {};
			\node at (5,10) [circ] {};
			\node at (6.25,10) [circ] {};
			\node at (7.5,10) [circ] {};
			\node at (8.75,10) [circ] {};
			\node at (10,10) [circ] {};
			\node at (11.25,10) [circ] {};
			\node at (12.5,10) [circ] {};
			\node at (1.25,11.25) [circ] {};
			\node at (2.5,11.25) [circ] {};
			\node at (3.75,11.25) [circ] {};
			\node at (5,11.25) [circ] {};
			\node at (6.25,11.25) [circ] {};
			\node at (7.5,11.25) [circ] {};
			\node at (8.75,11.25) [circ] {};
			\node at (10,11.25) [circ] {};
			\node at (11.25,11.25) [circ] {};
			\draw [ line width=1pt](3.75,8.75) to[short] (5,8.75);
			\draw [ line width=1pt](6.25,10) to[short] (6.25,8.75);
			\draw [line width=1pt, short] (5,8.75) -- (6.25,10);
			\draw [line width=1pt, short] (6.25,10) -- (7.5,8.75);
			\draw [line width=1pt, short] (6.25,10) -- (8.75,10);
			\draw [line width=1pt, short] (8.75,10) -- (12.5,10);
			\draw [line width=1pt, short] (7.5,10) -- (8.75,11.25);
			\draw [line width=1pt, short] (8.75,11.25) -- (11.25,11.25);
			\draw [line width=1pt, short] (6.25,10) -- (6.25,11.25);
			\draw [line width=1pt, short] (6.25,11.25) -- (5,10);
			\draw [line width=1pt, short] (5,11.25) -- (6.25,11.25);
			\draw [line width=1pt, short] (3.75,11.25) -- (5,11.25);
			\draw [line width=1pt, short] (2.5,11.25) -- (3.75,11.25);
			\draw [line width=1pt, short] (2.5,10) -- (3.75,11.25);
			\draw [line width=1pt, short] (3.75,11.25) -- (3.75,10);
			\node at (0,12.5) [circ] {};
			\node at (1.25,12.5) [circ] {};
			\node at (3.75,12.5) [circ] {};
			\node at (5,12.5) [circ] {};
			\node at (6.25,12.5) [circ] {};
			\node at (7.5,12.5) [circ] {};
			\node at (10,12.5) [circ] {};
			\node at (0,13.75) [circ] {};
			\node at (0,15) [circ] {};
			\node at (1.25,15) [circ] {};
			\node at (1.25,13.75) [circ] {};
			\node at (2.5,13.75) [circ] {};
			\node at (3.75,13.75) [circ] {};
			\node at (3.75,15) [circ] {};
			\node at (5,13.75) [circ] {};
			\node at (6.25,15) [circ] {};
			\node at (7.5,15) [circ] {};
			\node at (8.75,13.75) [circ] {};
			\node at (8.75,15) [circ] {};
			\node at (10,15) [circ] {};
			\draw [line width=1pt, short] (0,13.75) -- (1.25,13.75);
			\draw [line width=1pt, short] (0,12.5) -- (1.25,13.75);
			\draw [line width=1pt, short] (0,15) -- (1.25,13.75);
			\draw [line width=1pt, short] (1.25,13.75) -- (2.5,13.75);
			\draw [line width=1pt, short] (1.25,15) -- (2.5,13.75);
			\draw [line width=1pt, short] (1.25,12.5) -- (2.5,13.75);
			\draw [line width=1pt, short] (2.5,13.75) -- (3.75,13.75);
			\draw [line width=1pt, short] (3.75,15) -- (3.75,13.75);
			\draw [line width=1pt, short] (3.75,13.75) -- (5,12.5);
			\draw [line width=1pt, short] (3.75,12.5) -- (5,12.5);
			\draw [line width=1pt, short] (5,13.75) -- (5,12.5);
			\draw [line width=1pt, short] (5,12.5) -- (6.25,11.25);
			\draw [line width=1pt, short] (6.25,12.5) -- (6.25,11.25);
			\draw [line width=1pt, short] (6.25,12.5) -- (7.5,12.5);
			\draw [line width=1pt, short] (6.25,12.5) -- (6.25,15);
			\draw [line width=1pt, short] (6.25,15) -- (7.5,15);
			\draw [line width=1pt, short] (7.5,15) -- (8.75,13.75);
			\draw [line width=1pt, short] (7.5,15) -- (8.75,15);
			\draw [line width=1pt, short] (1.25,11.25) .. controls (2.75,12) and (2.5,11.75) .. (3.75,11.25);
			\draw [line width=1pt, short] (6.25,12.5) -- (7.5,11.25);
			\draw [line width=1pt, short] (7.5,15) .. controls (8.75,15.5) and (8.75,15.5) .. (10,15);
			\node at (11.25,15) [circ] {};
			\node at (12.5,15) [circ] {};
			\node at (12.5,13.75) [circ] {};
			\node at (13.75,15) [circ] {};
			\node at (6.25,16.25) [circ] {};
			\node at (7.5,16.25) [circ] {};
			\node at (8.75,16.25) [circ] {};
			\node at (10,16.25) [circ] {};
			\node at (11.25,16.25) [circ] {};
			\node at (12.5,16.25) [circ] {};
			\node at (13.75,16.25) [circ] {};
			\node at (15,16.25) [circ] {};
			\node at (16.25,16.25) [circ] {};
			\node at (16.25,17.5) [circ] {};
			\draw [line width=1pt, short] (6.25,16.25) -- (6.25,15);
			\draw [line width=1pt, short] (6.25,15) -- (7.5,16.25);
			\draw [line width=1pt, short] (10,16.25) -- (11.25,15);
			\draw [line width=1pt, short] (10,16.25) -- (11.25,16.25);
			\draw [line width=1pt, short] (11.25,16.25) -- (12.5,16.25);
			\draw [line width=1pt, short] (12.5,16.25) -- (13.75,16.25);
			\draw [line width=1pt, short] (15,16.25) -- (16.25,16.25);
			\draw [line width=1pt, short] (15,16.25) -- (16.25,17.5);
			\draw [line width=1pt, short] (12.5,16.25) -- (13.75,15);
			\draw [line width=1pt, short] (11.25,16.25) -- (12.5,15);
			\draw [line width=1pt, short] (11.25,16.25) -- (12.5,13.75);
			\node at (6.25,17.5) [circ] {};
			\node at (7.5,17.5) [circ] {};
			\node at (8.75,17.5) [circ] {};
			\node at (10,17.5) [circ] {};
			\node at (11.25,17.5) [circ] {};
			\node at (12.5,17.5) [circ] {};
			\draw [line width=1pt, short] (6.25,17.5) -- (6.25,16.25);
			\draw [line width=1pt, short] (7.5,17.5) -- (8.75,17.5);
			\draw [line width=1pt, short] (8.75,17.5) -- (8.75,16.25);
			\draw [line width=1pt, short] (8.75,17.5) -- (10,16.25);
			\draw [line width=1pt, short] (8.75,17.5) -- (10,17.5);
			\draw [line width=1pt, short] (10,17.5) -- (11.25,17.5);
			\draw [line width=1pt, short] (11.25,17.5) -- (12.5,17.5);
			\draw [line width=1pt, short] (6.25,16.25) -- (7.5,17.5);
			\node at (6.25,20) [circ] {};
			\node at (7.5,18.75) [circ] {};
			\node at (7.5,20) [circ] {};
			\node at (8.75,18.75) [circ] {};
			\node at (8.75,20) [circ] {};
			\node at (10,18.75) [circ] {};
			\node at (11.25,18.75) [circ] {};
			\node at (11.25,20) [circ] {};
			\node at (12.5,20) [circ] {};
			\node at (6.25,21.25) [circ] {};
			\node at (7.5,21.25) [circ] {};
			\node at (8.75,21.25) [circ] {};
			\node at (10,21.25) [circ] {};
			\node at (12.5,21.25) [circ] {};
			\node at (6.25,22.5) [circ] {};
			\node at (8.75,22.5) [circ] {};
			\node at (11.25,22.5) [circ] {};
			\node at (12.5,22.5) [circ] {};
			\node at (12.5,23.75) [circ] {};
			\node at (11.25,23.75) [circ] {};
			\draw [line width=1pt, short] (6.25,20) -- (6.25,17.5);
			\draw [line width=1pt, short] (6.25,21.25) -- (6.25,20);
			\draw [line width=1pt, short] (6.25,20) -- (7.5,21.25);
			\draw [line width=1pt, short] (6.25,20) -- (7.5,20);
			\draw [line width=1pt, short] (6.25,22.5) -- (7.5,21.25);
			\draw [line width=1pt, short] (7.5,21.25) -- (8.75,22.5);
			\draw [line width=1pt, short] (7.5,21.25) -- (8.75,21.25);
			\draw [line width=1pt, short] (8.75,22.5) -- (10,21.25);
			\draw [line width=1pt, short] (8.75,22.5) -- (11.25,22.5);
			\draw [line width=1pt, short] (11.25,23.75) -- (11.25,22.5);
			\draw [line width=1pt, short] (11.25,22.5) -- (12.5,23.75);
			\draw [line width=1pt, short] (11.25,22.5) -- (12.5,22.5);
			\draw [line width=1pt, short] (11.25,22.5) -- (12.5,21.25);
			\draw [line width=1pt, short] (7.5,18.75) -- (7.5,17.5);
			\draw [line width=1pt, short] (7.5,17.5) -- (8.75,18.75);
			\draw [line width=1pt, short] (8.75,17.5) -- (10,18.75);
			\draw [line width=1pt, short] (10,17.5) -- (11.25,18.75);
			\draw [line width=1pt, short] (11.25,20) -- (11.25,18.75);
			\draw [line width=1pt, short] (11.25,18.75) -- (12.5,20);
			\draw [line width=1pt, short] (10,16.25) .. controls (12.5,16.75) and (13.25,17.25) .. (15,16.25);
			\draw [line width=1pt, short] (6.25,20) .. controls (7.5,20.5) and (7.5,20.5) .. (8.75,20);
			\node at (5,16.25) [circ] {};
			\node at (2.5,16.25) [circ] {};
			\node at (1.25,16.25) [circ] {};
			\node at (0,16.25) [circ] {};
			\node at (-1.25,16.25) [circ] {};
			\node at (-2.5,16.25) [circ] {};
			\node at (-3.75,16.25) [circ] {};
			\node at (-5,16.25) [circ] {};
			\node at (5,17.5) [circ] {};
			\node at (3.75,17.5) [circ] {};
			\node at (0,17.5) [circ] {};
			\node at (-2.5,17.5) [circ] {};
			\node at (-3.75,17.5) [circ] {};
			\node at (-3.75,18.75) [circ] {};
			\node at (-2.5,18.75) [circ] {};
			\node at (-3.75,20) [circ] {};
			\node at (-2.5,20) [circ] {};
			\node at (-1.25,18.75) [circ] {};
			\node at (0,18.75) [circ] {};
			\node at (1.25,18.75) [circ] {};
			\node at (3.75,18.75) [circ] {};
			\node at (5,18.75) [circ] {};
			\node at (3.75,20) [circ] {};
			\node at (2.5,20) [circ] {};
			\node at (1.25,20) [circ] {};
			\node at (0,20) [circ] {};
			\draw [line width=1pt, short] (-3.75,17.5) -- (-2.5,17.5);
			\draw [line width=1pt, short] (-2.5,17.5) -- (-2.5,16.25);
			\draw [line width=1pt, short] (-5,16.25) -- (-2.5,17.5);
			\draw [line width=1pt, short] (-3.75,16.25) -- (-2.5,17.5);
			\draw [line width=1pt, short] (-2.5,17.5) -- (0,17.5);
			\draw [line width=1pt, short] (-1.25,16.25) -- (0,17.5);
			\draw [line width=1pt, short] (0,16.25) -- (2.5,16.25);
			\draw [line width=1pt, short] (1.25,16.25) -- (3.75,17.5);
			\draw [line width=1pt, short] (0,17.5) .. controls (1.5,18.25) and (2.25,18.5) .. (3.75,17.5);
			\draw [line width=1pt, short] (5,18.75) -- (3.75,17.5);
			\draw [line width=1pt, short] (5,17.5) -- (6.25,17.5);
			\draw [line width=1pt, short] (5,18.75) -- (6.25,17.5);
			\draw [line width=1pt, short] (3.75,18.75) -- (5,18.75);
			\draw [line width=1pt, short] (5,16.25) -- (6.25,15);
			\draw [line width=1pt, short] (3.75,20) -- (3.75,18.75);
			\draw [line width=1pt, short] (2.5,20) -- (3.75,18.75);
			\draw [line width=1pt, short] (1.25,18.75) -- (3.75,18.75);
			\draw [line width=1pt, short] (1.25,20) -- (1.25,18.75);
			\draw [line width=1pt, short] (0,20) -- (1.25,18.75);
			\draw [line width=1pt, short] (-1.25,18.75) -- (1.25,18.75);
			\draw [line width=1pt, short] (-2.5,18.75) -- (0,17.5);
			\draw [line width=1pt, short] (-2.5,20) -- (-2.5,18.75);
			\draw [line width=1pt, short] (-3.75,18.75) -- (-2.5,18.75);
			\node at (1.25,21.25) [circ] {};
			\node at (2.5,21.25) [circ] {};
			\node at (3.75,21.25) [circ] {};
			\node at (3.75,22.5) [circ] {};
			\node at (2.5,22.5) [circ] {};
			\node at (1.25,22.5) [circ] {};
			\draw [line width=1pt, short] (-3.75,20) -- (-2.5,18.75);
			\draw [line width=1pt, short] (1.25,22.5) -- (2.5,21.25);
			\draw [line width=1pt, short] (1.25,21.25) -- (2.5,21.25);
			\draw [line width=1pt, short] (2.5,22.5) -- (2.5,21.25);
			\draw [line width=1pt, short] (3.75,22.5) -- (3.75,21.25);
			\draw [line width=1pt, short] (2.5,21.25) -- (3.75,21.25);
			\draw [line width=1pt, short] (3.75,21.25) -- (5,18.75);
			\draw [line width=1pt, short] (8.75,11.25) -- (10,12.5);
			\node [font=\Large] at (6,9.5) {};
			\node [font=\Large] at (6.0,10.1) {0};
			\node [font=\Large] at (6.25,8.4) {32};
			\node [font=\Large] at (7.75,8.5) {2};
			\node [font=\Large] at (5,8.5) {64};
			\node [font=\Large] at (3.5,8.5) {65};
			\node [font=\Large] at (8.5,11.5) {};
			\node [font=\Large] at (2.5,9.5) {};
			\node [font=\Large] at (3.5,9.75) {86};
			\node [font=\Large] at (2.25,9.75) {66};
			\node [font=\Large] at (4.75,9.75) {48};
			\node [font=\Large] at (7.5,9.6) {8};
			\node [font=\Large] at (8.5,9.70) {40};
			\node [font=\Large] at (9.75,9.70) {44};
			\node [font=\Large] at (11,9.70) {45};
			\node [font=\Large] at (8.5,11.5) {};
			\node [font=\Large] at (8.5,11.5) {};
			\node [font=\Large] at (12.5,9.70) {37};
			\node [font=\Large] at (9,10.9) {72};
			\node [font=\Large] at (10,10.9) {74};
			\node [font=\Large] at (11.25,10.9) {10};
			\node [font=\Large] at (10.25,12.25) {104};
			\node [font=\Large] at (1.25,11) {18};
			\node [font=\Large] at (2.25,11) {90};
			\node [font=\Large] at (4,11.5) {82};
			\node [font=\Large] at (5,10.9) {80};
			\node [font=\Large] at (6.55,11.25) {16};
			\node [font=\Large] at (5.9,12.5) {17};
			\node [font=\Large] at (7.5,10.9) {25};
			\node [font=\Large] at (7.75,12.5) {1};
			\node [font=\Large] at (8.75,13.5) {93};
			\node [font=\Large] at (6,14.75) {21};
			\node [font=\Large] at (7.2,14.75) {85};
			\node [font=\Large] at (8.75,14.75) {81};
			\node [font=\Large] at (9.75,14.75) {87};
			\node [font=\Large] at (11,14.75) {75};
			\node [font=\Large] at (12.25,13.5) {109};
			\node [font=\Large] at (12.5,14.75) {107};
			\node [font=\Large] at (13.75,14.75) {108};
			\node [font=\Large] at (16,17.75) {73};
			\node [font=\Large] at (16.25,16.5) {76};
			\node [font=\Large] at (15,15.9) {77};
			\node [font=\Large] at (9.70,15.95) {79};
			\node [font=\Large] at (10.9,15.9) {111};
			\node [font=\Large] at (12.5,16.5) {110};
			\node [font=\Large] at (13.75,15.9) {78};
			\node [font=\Large] at (8.5,15.9) {11};
			\node [font=\Large] at (7.5,16) {};
			\node [font=\Large] at (6.7,16.25) {23};
			\node [font=\Large] at (4.75,16) {20};
			\node [font=\Large] at (12.75,17.5) {4};
			\node [font=\Large] at (11.25,17.25) {};
			\node [font=\Large] at (10,17.25) {};
			\node [font=\Large] at (8.55,17.75) {15};
			\node [font=\Large] at (7.5,17) {31};
			\node [font=\Large] at (6.65,17.5) {55};
			\node [font=\Large] at (7,18.75) {29};
			\node [font=\Large] at (8.5,19) {30};
			\node [font=\Large] at (10,19.25) {14};
			\node [font=\Large] at (11.25,18.25) {71};
			\node [font=\Large] at (12.75,20) {  67};
			\node [font=\Large] at (11.5,20.25) {69};
			\node [font=\Large] at (9.25,20) {50};
			\node [font=\Large] at (7.75,19.75) {52};
			\node [font=\Large] at (6.55,19.75) {54};
			\node [font=\Large] at (5.75,21.25) {62};
			\node [font=\Large] at (7.65,20.9) {38};
			\node [font=\Large] at (8.75,20.75) {39};
			\node [font=\Large] at (10,21) {118};
			\node [font=\Large] at (13,21.25) {116};
			\node [font=\Large] at (13,22.5) {101};
			\node [font=\Large] at (12.75,23.5) {68};
			\node [font=\Large] at (10.75,23.75) {96};
			\node [font=\Large] at (10.75,22.75) {100};
			\node [font=\Large] at (8.75,22.75) {102};
			\node [font=\Large] at (6.25,22.75) {34};
			\node [font=\Large] at (5,17.2) {53};
			\node [font=\Large] at (3.75,22.75) {119};
			\node [font=\Large] at (4.25,21.25) {  115};
			\node [font=\Large] at (2.5,20.75) {  113};
			\node [font=\Large] at (2.5,23) {117};
			\node [font=\Large] at (0.75,22.25) {112};
			\node [font=\Large] at (0.75,21.25) {49};
			\node [font=\Large] at (-4,19.6) {114};
			\node [font=\Large] at (-2.25,20.3) {   106};
			\node [font=\Large] at (-4.25,18.70) {58};
			\node [font=\Large] at (-2.75,18.4) {122};
			\node [font=\Large] at (-1.25,19) {105};
			\node [font=\Large] at (0,19) {97};
			\node [font=\Large] at (1.25,18.4) {99};
			\node [font=\Large] at (-0.5,20) {98};
			\node [font=\Large] at (1,20) {};
			\node [font=\Large] at (0.75,20) {103};
			\node [font=\Large] at (2.25,19.75) {3};
			\node [font=\Large] at (3.4,20) {33};
			\node [font=\Large] at (5.4,18.75) {  51};
			\node [font=\Large] at (3.4,18.45) {35};
			\node [font=\Large] at (3.75,17.9) {59};
			\node [font=\Large] at (0.25,17.25) {123};
			\node [font=\Large] at (-2.2,17.20) {91};
			\node [font=\Large] at (-4.25,17.5) {95};
			\node [font=\Large] at (-5,16.6) {89};
			\node [font=\Large] at (-4,15.9) {83};
			\node [font=\Large] at (-2.85,16.25) {27};
			\node [font=\Large] at (-1.25,15.9) {121};
			\node [font=\Large] at (0,15.9) {47};
			\node [font=\Large] at (1.2,16.55) {63};
			\node [font=\Large] at (2.95,16.25) {127};
			\node [font=\Large] at (-0.5,15) {125};
			\node [font=\Large] at (-0.5,13.75) {60};
			\node [font=\Large] at (-0.5,12.5) {126};
			\node [font=\Large] at (1,12.25) {94};
			\node [font=\Large] at (1.5,13.4) {124};
			\node [font=\Large] at (2.65,14) {92};
			\node [font=\Large] at (1.65,15) {84};
			\node [font=\Large] at (3.25,15) {120};
			\node [font=\Large] at (4.1,14) {88};
			\node [font=\Large] at (5.25,13.5) {26};
			\node [font=\Large] at (4.75,12.15) {24};
			\node [font=\Large] at (3.25,12.5) {28};
			\node [font=\Large] at (10,17.2) {7};
			\node [font=\Large] at (11.25,17.2) {6};
			\node [font=\Large] at (7.6,15.9) {5};
			\node at (8.75,8.75) [circ] {};
			\node at (10,7.5) [circ] {};
			\node at (10,8.75) [circ] {};
			\node at (11.25,8.75) [circ] {};
			\node at (11.25,7.5) [circ] {};
			\node at (12.5,8.75) [circ] {};
			\node at (12.5,7.5) [circ] {};
			\node at (13.75,10) [circ] {};
			\draw [line width=1pt, short] (7.5,10) -- (8.75,8.75);
			\draw [line width=1pt, short] (8.75,10) -- (10,8.75);
			\draw [line width=1pt, short] (8.75,10) -- (10,7.5);
			\draw [line width=1pt, short] (10,10) -- (11.25,8.75);
			\draw [line width=1pt, short] (10,10) -- (11.25,7.5);
			\draw [line width=1pt, short] (11.25,10) -- (12.5,7.5);
			\draw [line width=1pt, short] (11.25,10) -- (12.5,8.75);
			\draw [line width=1pt, short] (11.25,10) .. controls (12.5,10.5) and (12.5,10.5) .. (13.75,10);
			\node [font=\Large] at (8.75,8.25) {9};
			\node [font=\Large] at (9.5,7.5) {56};
			\node [font=\Large] at (10,8.5) {42};
			\node [font=\Large] at (10.75,7.5) {36};
			\node [font=\Large] at (11.25,8.5) {12};
			\node [font=\Large] at (12,7.5) {13};
			\node [font=\Large] at (12.75,8.5) {61};
			\node [font=\Large] at (13.75,9.5) {41};
			\node at (3.75,16.25) [circ] {};
			\draw [ line width=1pt](3.75,17.5) to[short] (3.75,16.25);
			\node [font=\Large] at (4.1,16.5) {57};
			\node at (5,20) [circ] {};
			\node at (12.5,18.75) [circ] {};
			\node at (13.75,17.5) [circ] {};
			\node at (15,15) [circ] {};
			\draw [line width=1pt, short] (5,20) -- (5,18.75);
			\draw [line width=1pt, short] (11.25,18.75) -- (12.5,18.75);
			\draw [line width=1pt, short] (12.5,16.25) -- (15,15);
			\draw [line width=1pt, short] (11.25,17.5) .. controls (12.75,18.25) and (12.75,18.25) .. (13.75,17.5);
			\node [font=\Large] at (5.25,20.5) {19};
			\node [font=\Large] at (13,19) {70};
			\node [font=\Large] at (14,17.75) {22};
			\node [font=\Large] at (15.5,15) {46};
			\node at (2,17.5) [circ] {};
			\draw [line width=1pt, short] (2,17.5) -- (3.75,17.5);
			\node [font=\Large] at (1.5,17.5) {43};
		\end{circuitikz}
	}%

		\caption{$T_1:$ First CIST in $Q_7$ with diameter 15}
	\end{figure}

	\vspace*{-0.5 cm}	
	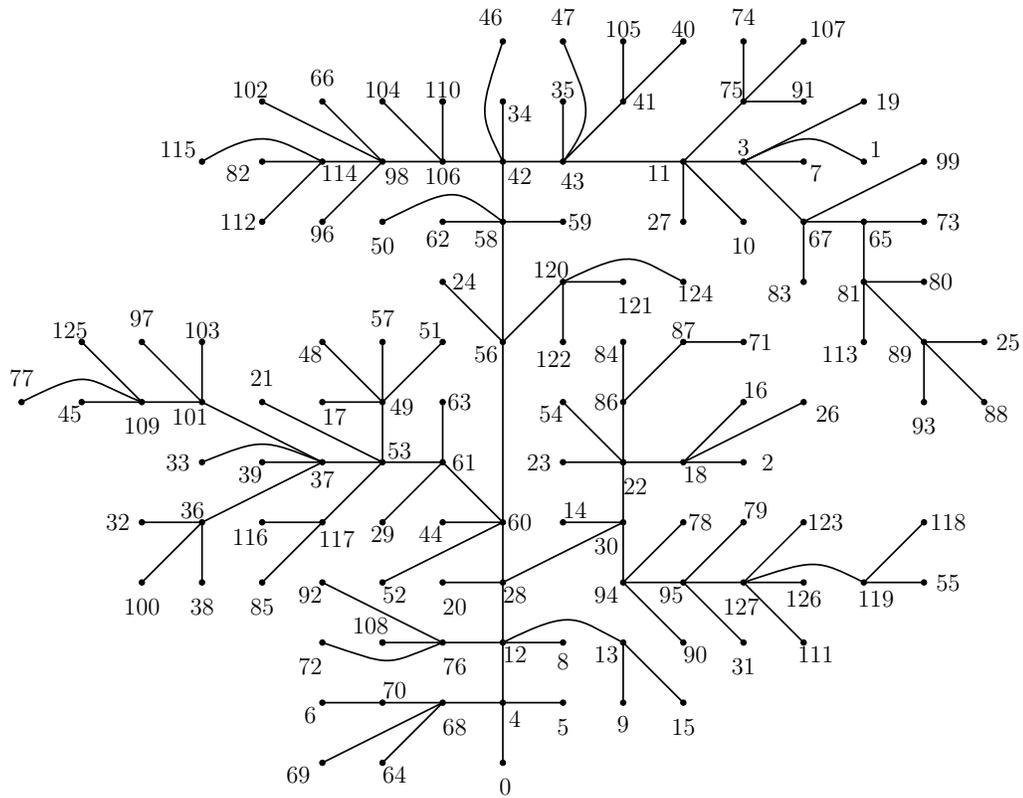
\begin{figure}[!ht]
		\centering
		\resizebox{0.9\textwidth}{!}{%
			\begin{circuitikz}
				\tikzstyle{every node}=[font=\Large]
				\node at (23.75,20) [circ] {};
				\node at (21.25,20) [circ] {};
				\node at (20,20) [circ] {};
				\node at (23.75,21.25) [circ] {};
				\node at (22.5,21.25) [circ] {};
				\node at (21.25,21.25) [circ] {};
				\node at (25,21.25) [circ] {};
				\node at (26.25,21.25) [circ] {};
				\node at (27.5,21.25) [circ] {};
				\draw [line width=1pt, short] (23.75,21.25) -- (23.75,20);
				\draw [line width=1pt, short] (21.25,21.25) -- (25,21.25);
				\draw [line width=1pt, short] (20,20) -- (22.5,21.25);
				\draw [line width=1pt, short] (21.25,20) -- (22.5,21.25);
				\node at (23.75,22.5) [circ] {};
				\node at (22.5,22.5) [circ] {};
				\node at (21.25,22.5) [circ] {};
				\node at (20,22.5) [circ] {};
				\node at (21.25,22.5) [circ] {};
				\node at (25,22.5) [circ] {};
				\node at (26.25,22.5) [circ] {};
				\node at (27.5,22.5) [circ] {};
				\node at (28.75,22.5) [circ] {};
				\node at (30,22.5) [circ] {};
				\draw [line width=1pt, short] (21.25,22.5) -- (22.5,22.5);
				\draw [line width=1pt, short] (22.5,22.5) -- (23.75,22.5);
				\draw [line width=1pt, short] (23.75,22.5) -- (23.75,21.25);
				\draw [line width=1pt, short] (23.75,22.5) -- (25,22.5);
				\draw [line width=1pt, short] (26.25,22.5) -- (26.25,21.25);
				\draw [line width=1pt, short] (26.25,22.5) -- (27.5,21.25);
				\node at (23.75,23.75) [circ] {};
				\node at (22.5,23.75) [circ] {};
				\node at (21.25,23.75) [circ] {};
				\node at (20,23.75) [circ] {};
				\node at (18.75,23.75) [circ] {};
				\node at (17.5,23.75) [circ] {};
				\node at (16.25,23.75) [circ] {};
				\node at (26.25,23.75) [circ] {};
				\node at (27.5,23.75) [circ] {};
				\node at (28.75,23.75) [circ] {};
				\node at (30,23.75) [circ] {};
				\node at (31.25,23.75) [circ] {};
				\node at (23.75,25) [circ] {};
				\node at (22.5,25) [circ] {};
				\node at (21.25,25) [circ] {};
				\node at (20,25) [circ] {};
				\node at (18.75,25) [circ] {};
				\node at (17.5,25) [circ] {};
				\node at (16.25,25) [circ] {};
				\node at (25,25) [circ] {};
				\node at (26.25,25) [circ] {};
				\node at (27.5,25) [circ] {};
				\node at (28.75,25) [circ] {};
				\node at (30,25) [circ] {};
				\draw [line width=1pt, short] (23.75,22.5) -- (23.75,25);
				\draw [line width=1pt, short] (23.75,23.75) -- (22.5,23.75);
				\draw [line width=1pt, short] (23.75,25) -- (21.25,23.75);
				\draw [line width=1pt, short] (22.5,22.5) -- (20,23.75);
				\draw [line width=1pt, short] (20,25) -- (18.75,23.75);
				\draw [line width=1pt, short] (18.75,25) -- (20,25);
				\draw [line width=1pt, short] (16.25,25) -- (17.5,25);
				\draw [line width=1pt, short] (16.25,23.75) -- (17.5,25);
				\draw [line width=1pt, short] (17.5,25) -- (17.5,23.75);
				\draw [line width=1pt, short] (22.5,25) -- (23.75,25);
				\draw [line width=1pt, short] (26.25,25) -- (26.25,23.75);
				\draw [line width=1pt, short] (25,25) -- (26.25,25);
				\draw [line width=1pt, short] (23.75,23.75) -- (26.25,25);
				\draw [line width=1pt, short] (26.25,23.75) -- (30,23.75);
				\draw [line width=1pt, short] (26.25,23.75) -- (27.5,25);
				\draw [line width=1pt, short] (26.25,23.75) -- (27.5,22.5);
				\draw [line width=1pt, short] (27.5,23.75) -- (28.75,25);
				\draw [line width=1pt, short] (27.5,23.75) -- (28.75,22.5);
				\draw [line width=1pt, short] (28.75,23.75) -- (30,25);
				\draw [line width=1pt, short] (28.75,23.75) -- (30,22.5);
				\node at (32.5,23.75) [circ] {};
				\node at (32.5,25) [circ] {};
				\draw [line width=1pt, short] (31.25,23.75) -- (32.5,23.75);
				\draw [line width=1pt, short] (32.5,25) -- (31.25,23.75);
				\node at (23.75,28.75) [circ] {};
				\node at (23.75,31.25) [circ] {};
				\node at (23.75,32.5) [circ] {};
				\node at (23.75,33.75) [circ] {};
				\node at (23.75,35) [circ] {};
				\node at (22.5,26.25) [circ] {};
				\node at (21.25,26.25) [circ] {};
				\node at (20,26.25) [circ] {};
				\node at (18.75,26.25) [circ] {};
				\node at (17.5,26.25) [circ] {};
				\node at (22.5,27.5) [circ] {};
				\node at (21.25,27.5) [circ] {};
				\node at (20,27.5) [circ] {};
				\node at (18.75,27.5) [circ] {};
				\node at (17.5,27.5) [circ] {};
				\node at (16.25,27.5) [circ] {};
				\node at (15,27.5) [circ] {};
				\node at (13.75,27.5) [circ] {};
				\node at (22.5,28.75) [circ] {};
				\node at (21.25,28.75) [circ] {};
				\node at (20,28.75) [circ] {};
				\node at (17.5,28.75) [circ] {};
				\node at (16.25,28.75) [circ] {};
				\node at (15,28.75) [circ] {};
				\node at (22.5,30) [circ] {};
				\node at (22.5,31.25) [circ] {};
				\node at (21.25,31.25) [circ] {};
				\node at (20,31.25) [circ] {};
				\node at (18.75,31.25) [circ] {};
				\node at (22.5,32.5) [circ] {};
				\node at (21.25,32.5) [circ] {};
				\node at (20,32.5) [circ] {};
				\node at (18.75,32.5) [circ] {};
				\node at (17.5,32.5) [circ] {};
				\draw [line width=1pt, short] (23.75,33.75) -- (23.75,25);
				\draw [line width=1pt, short] (23.75,25) -- (22.5,26.25);
				\draw [line width=1pt, short] (22.5,26.25) -- (21.25,25);
				\draw [line width=1pt, short] (22.5,26.25) -- (21.25,26.25);
				\draw [line width=1pt, short] (22.5,26.25) -- (22.5,27.5);
				\draw [line width=1pt, short] (21.25,26.25) -- (21.25,27.5);
				\draw [line width=1pt, short] (21.25,26.25) -- (18.75,26.25);
				\draw [line width=1pt, short] (21.25,27.5) -- (20,27.5);
				\draw [line width=1pt, short] (21.25,26.25) -- (18.75,27.5);
				\draw [line width=1pt, short] (20,26.25) -- (17.5,27.5);
				\draw [line width=1pt, short] (17.5,27.5) -- (15,27.5);
				\draw [line width=1pt, short] (17.5,28.75) -- (17.5,27.5);
				\draw [line width=1pt, short] (16.25,28.75) -- (17.5,27.5);
				\draw [line width=1pt, short] (15,28.75) -- (16.25,27.5);
				\draw [line width=1pt, short] (21.25,27.5) -- (22.5,28.75);
				\draw [line width=1pt, short] (21.25,28.75) -- (21.25,27.5);
				\draw [line width=1pt, short] (20,28.75) -- (21.25,27.5);
				\draw [line width=1pt, short] (22.5,30) -- (23.75,28.75);
				\draw [line width=1pt, short] (22.5,31.25) -- (23.75,31.25);
				\draw [line width=1pt, short] (18.75,32.5) -- (23.75,32.5);
				\node at (22.5,33.75) [circ] {};
				\node at (21.25,33.75) [circ] {};
				\node at (20,33.75) [circ] {};
				\node at (18.75,33.75) [circ] {};
				\node at (22.5,33.75) [circ] {};
				\draw [line width=1pt, short] (22.5,33.75) -- (22.5,32.5);
				\draw [line width=1pt, short] (21.25,33.75) -- (22.5,32.5);
				\draw [line width=1pt, short] (20,33.75) -- (21.25,32.5);
				\draw [line width=1pt, short] (18.75,33.75) -- (21.25,32.5);
				\draw [line width=1pt, short] (20,32.5) -- (18.75,31.25);
				\draw [line width=1pt, short] (21.25,32.5) -- (20,31.25);
				\draw [line width=1pt, short] (17.5,25) -- (20,26.25);
				\draw [line width=1pt, short] (20,25) -- (21.25,26.25);
				\draw [line width=1pt, short] (20,22.5) .. controls (21.5,22) and (21.25,22) .. (22.5,22.5);
				\draw [line width=1pt, short] (17.5,26.25) .. controls (18.75,26.75) and (18.75,26.75) .. (20,26.25);
				\draw [line width=1pt, short] (13.75,27.5) .. controls (15.25,28.25) and (15,28) .. (16.25,27.5);
				\draw [line width=1pt, short] (21.25,31.25) .. controls (22.75,32) and (22.75,32) .. (23.75,31.25);
				\draw [line width=1pt, short] (17.5,32.5) .. controls (18.75,33.25) and (19,33) .. (20,32.5);
				\draw [line width=1pt, short] (23.75,35) .. controls (23.25,33.5) and (23.25,33.5) .. (23.75,32.5);
				\node at (25,26.25) [circ] {};
				\node at (26.25,26.25) [circ] {};
				\node at (27.5,26.25) [circ] {};
				\node at (28.75,26.25) [circ] {};
				\node at (25,27.5) [circ] {};
				\node at (26.25,27.5) [circ] {};
				\node at (28.75,27.5) [circ] {};
				\node at (30,27.5) [circ] {};
				\node at (32.5,27.5) [circ] {};
				\node at (33.75,27.5) [circ] {};
				\node at (25,28.75) [circ] {};
				\node at (26.25,28.75) [circ] {};
				\node at (27.5,28.75) [circ] {};
				\node at (28.75,28.75) [circ] {};
				\node at (32.5,28.75) [circ] {};
				\node at (33.75,28.75) [circ] {};
				\node at (25,30) [circ] {};
				\node at (26.25,30) [circ] {};
				\node at (27.5,30) [circ] {};
				\node at (30,30) [circ] {};
				\node at (31.25,30) [circ] {};
				\node at (32.5,30) [circ] {};
				\node at (25,31.25) [circ] {};
				\node at (27.5,31.25) [circ] {};
				\node at (28.75,31.25) [circ] {};
				\node at (30,31.25) [circ] {};
				\node at (31.25,31.25) [circ] {};
				\node at (32.5,31.25) [circ] {};
				\node at (25,32.5) [circ] {};
				\node at (27.5,32.5) [circ] {};
				\node at (28.75,32.5) [circ] {};
				\node at (30,32.5) [circ] {};
				\node at (31.25,32.5) [circ] {};
				\node at (32.5,32.5) [circ] {};
				\node at (25,33.75) [circ] {};
				\node at (26.25,33.75) [circ] {};
				\node at (28.75,33.75) [circ] {};
				\node at (30,33.75) [circ] {};
				\node at (31.25,33.75) [circ] {};
				\node at (25,35) [circ] {};
				\node at (26.25,35) [circ] {};
				\node at (27.5,35) [circ] {};
				\node at (28.75,35) [circ] {};
				\node at (30,35) [circ] {};
				\node at (26.25,26.25) [circ] {};
				\draw [line width=1pt, short] (26.25,26.25) -- (26.25,25);
				\draw [line width=1pt, short] (25,26.25) -- (28.75,26.25);
				\draw [line width=1pt, short] (26.25,28.75) -- (26.25,26.25);
				\draw [line width=1pt, short] (25,27.5) -- (26.25,26.25);
				\draw [line width=1pt, short] (27.5,26.25) -- (28.75,27.5);
				\draw [line width=1pt, short] (27.5,26.25) -- (30,27.5);
				\draw [line width=1pt, short] (26.25,27.5) -- (27.5,28.75);
				\draw [line width=1pt, short] (27.5,28.75) -- (28.75,28.75);
				\draw [line width=1pt, short] (23.75,28.75) -- (25,30);
				\draw [line width=1pt, short] (25,30) -- (25,28.75);
				\draw [line width=1pt, short] (25,30) -- (26.25,30);
				\draw [line width=1pt, short] (28.75,23.75) .. controls (30.25,24.25) and (30.25,24.25) .. (31.25,23.75);
				\draw [line width=1pt, short] (25,30) .. controls (26.5,30.75) and (26.5,30.5) .. (27.5,30);
				\draw [line width=1pt, short] (23.75,22.5) .. controls (25.25,23.25) and (25.25,23) .. (26.25,22.5);
				\draw [line width=1pt, short] (25,35) .. controls (25.5,33.75) and (25.75,33.5) .. (25,32.5);
				\draw [line width=1pt, short] (28.75,32.5) .. controls (30.25,33.25) and (30.25,33) .. (31.25,32.5);
				\draw [line width=1pt, short] (23.75,31.25) -- (25,31.25);
				\draw [line width=1pt, short] (23.75,32.5) -- (30,32.5);
				\draw [line width=1pt, short] (25,33.75) -- (25,32.5);
				\draw [line width=1pt, short] (25,32.5) -- (26.25,33.75);
				\draw [line width=1pt, short] (26.25,35) -- (26.25,33.75);
				\draw [line width=1pt, short] (26.25,33.75) -- (27.5,35);
				\draw [line width=1pt, short] (27.5,32.5) -- (28.75,33.75);
				\draw [line width=1pt, short] (28.75,35) -- (28.75,33.75);
				\draw [line width=1pt, short] (28.75,33.75) -- (30,35);
				\draw [line width=1pt, short] (28.75,33.75) -- (30,33.75);
				\draw [line width=1pt, short] (28.75,32.5) -- (31.25,33.75);
				\draw [line width=1pt, short] (27.5,32.5) -- (27.5,31.25);
				\draw [line width=1pt, short] (27.5,32.5) -- (28.75,31.25);
				\draw [line width=1pt, short] (28.75,32.5) -- (30,31.25);
				\draw [line width=1pt, short] (30,31.25) -- (32.5,32.5);
				\draw [line width=1pt, short] (30,31.25) -- (32.5,31.25);
				\draw [line width=1pt, short] (30,31.25) -- (30,30);
				\draw [line width=1pt, short] (31.25,31.25) -- (31.25,30);
				\draw [line width=1pt, short] (31.25,30) -- (32.5,30);
				\draw [line width=1pt, short] (31.25,30) -- (32.5,28.75);
				\draw [line width=1pt, short] (32.5,28.75) -- (33.75,28.75);
				\draw [line width=1pt, short] (32.5,28.75) -- (32.5,27.5);
				\draw [line width=1pt, short] (32.5,28.75) -- (33.75,27.5);
				\node [font=\Large] at (23.8,19.5) {0};
				\node [font=\Large] at (21.5,19.75) {64};
				\node [font=\Large] at (19.5,19.75) {69};
				\node [font=\Large] at (21.5,21.5) {70};
				\node [font=\Large] at (22.75,20.75) {68};
				\node [font=\Large] at (24,20.9) {4};
				\node [font=\Large] at (25,20.75) {5};
				\node [font=\Large] at (26.25,20.8) {9};
				\node [font=\Large] at (27.5,20.75) {15};
				\node [font=\Large] at (24,22.25) {12};
				\node [font=\Large] at (25,22.1) {8};
				\node [font=\Large] at (25.9,22.25) {13};
				\node [font=\Large] at (27.75,22.25) {90};
				\node [font=\Large] at (28.75,22) {31};
				\node [font=\Large] at (30.25,22.25) {   111};
				\node [font=\Large] at (22.75,22) {76};
				\node [font=\Large] at (21,22.85) {108};
				\node [font=\Large] at (19.75,22) {72};
				\node [font=\Large] at (19.75,23.5) {92};
				\node [font=\Large] at (16.25,23.25) {100};
				\node [font=\Large] at (17.5,23.25) {38};
				\node [font=\Large] at (15.75,25) {32};
				\node [font=\Large] at (18.5,24.70) {116};
				\node [font=\Large] at (17.3,25.25) {36};
				\node [font=\Large] at (17,26.25) {33};
				\node [font=\Large] at (18.5,26) {39};
				\node [font=\Large] at (20.30,24.65) {117};
				\node [font=\Large] at (18.75,23.25) {85};
				\node [font=\Large] at (24,23.5) {28};
				\node [font=\Large] at (22.75,23.25) {20};
				\node [font=\Large] at (21.5,23.5) {52};
				\node [font=\Large] at (21.25,24.75) {29};
				\node [font=\Large] at (22.25,24.75) {44};
				\node [font=\Large] at (24.1,25) {  60};
				\node [font=\Large] at (22.95,26.25) {  61};
				\node [font=\Large] at (22.85,27.5) {  63};
				\node [font=\Large] at (21.6,26.5) {53};
				\node [font=\Large] at (21.65,27.4) {49};
				\node [font=\Large] at (20.25,27.20) {17};
				\node [font=\Large] at (19.75,28.5) {48};
				\node [font=\Large] at (21.25,29.25) {57};
				\node [font=\Large] at (22.25,29) {51};
				\node [font=\Large] at (18.75,28) {21};
				\node [font=\Large] at (17.25,27.2) {101};
				\node [font=\Large] at (17.5,29) {103};
				\node [font=\Large] at (16.25,29.25) {97};
				\node [font=\Large] at (14.75,29) {125};
				\node [font=\Large] at (13.75,28) {77};
				\node [font=\Large] at (14.75,27.25) {45};
				\node [font=\Large] at (16.25,27) {109};
				\node [font=\Large] at (20,25.9) {     37};
				\node [font=\Large] at (23.4,28.5) {56};
				\node [font=\Large] at (22.95,30) {24};
				\node [font=\Large] at (23.4,30.9) {58};
				\node [font=\Large] at (22.4,30.9) {62};
				\node [font=\Large] at (21.25,30.75) {50};
				\node [font=\Large] at (20,31) {      96};
				\node [font=\Large] at (18.25,31.25) {112};
				\node [font=\Large] at (17,32.75) {115};
				\node [font=\Large] at (20.35,32.25) {  114};
				\node [font=\Large] at (21.55,32.20) {98};
				\node [font=\Large] at (18.5,34) {102};
				\node [font=\Large] at (20,34.25) {66};
				\node [font=\Large] at (21.25,34) {104};
				\node [font=\Large] at (22.5,34) {110};
				\node [font=\Large] at (22.5,32.20) {106};
				\node [font=\Large] at (24.1,32.20) {42};
				\node [font=\Large] at (25.2,32.1) {43};
				\node [font=\Large] at (24.1,33.5) {34};
				\node [font=\Large] at (25,34) {35};
				\node [font=\Large] at (23.5,35.5) {46};
				\node [font=\Large] at (25,35.5) {47};
				\node [font=\Large] at (26.25,35.25) {105};
				\node [font=\Large] at (27.5,35.25) {40};
				\node [font=\Large] at (26.7,33.75) {   41};
				\node [font=\Large] at (27,32.25) {   11};
				\node [font=\Large] at (27,31.25) {27};
				\node [font=\Large] at (28.75,30.75) {10};
				\node [font=\Large] at (25.35,31.25) {  59};
				\node [font=\Large] at (28.75,35.5) {74};
				\node [font=\Large] at (30.5,35.25) {107};
				\node [font=\Large] at (30,34) {91};
				\node [font=\Large] at (31.75,33.75) {19};
				\node [font=\Large] at (31.5,32.75) {1};
				\node [font=\Large] at (30.25,32.25) {7};
				\node [font=\Large] at (28.75,32.8) {3};
				\node [font=\Large] at (28.5,34) {75};
				\node [font=\Large] at (33,32.5) {99};
				\node [font=\Large] at (33,31.25) {73};
				\node [font=\Large] at (31.6,30.9) {  65};
				\node [font=\Large] at (30.35,30.9) {  67};
				\node [font=\Large] at (29.5,29.75) {83};
				\node [font=\Large] at (30.95,29.75) {81};
				\node [font=\Large] at (32.85,30) {80};
				\node [font=\Large] at (32,28.5) {89};
				\node [font=\Large] at (34.25,28.75) {25};
				\node [font=\Large] at (34,27.25) {88};
				\node [font=\Large] at (32.5,27) {93};
				\node [font=\Large] at (30.5,27.25) {26};
				\node [font=\Large] at (29.25,26.25) {2};
				\node [font=\Large] at (29,27.75) {16};
				\node [font=\Large] at (27.75,26) {18};
				\node [font=\Large] at (26.5,25.75) {22};
				\node [font=\Large] at (25.9,24.5) {   30};
				\node [font=\Large] at (24.5,26.25) {23};
				\node [font=\Large] at (25.25,25.25) {14};
				\node [font=\Large] at (24.75,27.25) {54};
				\node [font=\Large] at (25.90,27.5) {86};
				\node [font=\Large] at (25.90,28.5) {84};
				\node [font=\Large] at (24.8,28.4) {122};
				\node [font=\Large] at (26.5,29.5) {121};
				\node [font=\Large] at (24.75,30.25) {120};
				\node [font=\Large] at (27.75,29.75) {    124};
				\node [font=\Large] at (27.5,29) {87};
				\node [font=\Large] at (29.1,28.75) {71};
				\node [font=\Large] at (27.85,25) {78};
				\node [font=\Large] at (29,25.25) {79};
				\node [font=\Large] at (30.45,25) {123};
				\node [font=\Large] at (33,25) {118};
				\node [font=\Large] at (33,23.75) {55};
				\node [font=\Large] at (31.5,23.4) {   119};
				\node [font=\Large] at (30,23.4) {126};
				\node [font=\Large] at (28.7,23.25) {127};
				\node [font=\Large] at (27.25,23.5) {95};
				\node [font=\Large] at (25.9,23.5) {   94};
				\node at (20,21.25) [circ] {};
				\node at (31.25,28.75) [circ] {};
				\draw [line width=1pt, short] (20,21.25) -- (21.25,21.25);
				\draw [line width=1pt, short] (31.25,30) -- (31.25,28.75);
				\node [font=\Large] at (30.75,28.5) {113};
				\node [font=\Large] at (19.75,21) {6};
				\node [font=\Large] at (18.25,32.25) {82};
			\end{circuitikz}
		}%

		\caption {$T_2$: Second CIST in Q7 with diameter 18}
	\end{figure}

	\vspace*{2 cm}	
	
		
		\begin{figure}[!ht]
			\centering
			\resizebox{1\textwidth}{!}{%
				\begin{circuitikz}
					\tikzstyle{every node}=[font=\Large]
					\node at (1.25,4.5) [circ] {};
					\node at (2.5,4.5) [circ] {};
					\node at (5,4.5) [circ] {};
					\node at (6.25,4.5) [circ] {};
					\node at (7.5,4.5) [circ] {};
					\node at (8.75,4.5) [circ] {};
					\node at (10,4.5) [circ] {};
					\node at (11.25,4.5) [circ] {};
					\node at (12.5,4.5) [circ] {};
					\node at (13.75,4.5) [circ] {};
					\node at (8.75,3.25) [circ] {};
					\node at (0,5.75) [circ] {};
					\node at (1.25,5.75) [circ] {};
					\node at (2.5,5.75) [circ] {};
					\node at (3.75,5.75) [circ] {};
					\node at (5,5.75) [circ] {};
					\node at (6.25,5.75) [circ] {};
					\node at (7.5,5.75) [circ] {};
					\node at (8.75,5.75) [circ] {};
					\node at (10,5.75) [circ] {};
					\node at (12.5,5.75) [circ] {};
					\node at (13.75,5.75) [circ] {};
					\node at (15,5.75) [circ] {};
					\node at (16.25,5.75) [circ] {};
					\node at (0,7) [circ] {};
					\node at (-1.25,7) [circ] {};
					\node at (-2.5,7) [circ] {};
					\node at (-3.75,7) [circ] {};
					\node at (1.25,7) [circ] {};
					\node at (2.5,7) [circ] {};
					\node at (3.75,7) [circ] {};
					\node at (5,7) [circ] {};
					\node at (6.25,7) [circ] {};
					\node at (7.5,7) [circ] {};
					\node at (8.75,7) [circ] {};
					\node at (10,7) [circ] {};
					\node at (11.25,7) [circ] {};
					\node at (12.5,7) [circ] {};
					\node at (13.75,7) [circ] {};
					\node at (15,7) [circ] {};
					\node at (-5,8.25) [circ] {};
					\node at (-3.75,8.25) [circ] {};
					\node at (-2.5,8.25) [circ] {};
					\node at (-1.25,8.25) [circ] {};
					\node at (0,8.25) [circ] {};
					\node at (2.5,8.25) [circ] {};
					\node at (3.75,8.25) [circ] {};
					\node at (5,8.25) [circ] {};
					\node at (6.25,8.25) [circ] {};
					\node at (7.5,8.25) [circ] {};
					\node at (8.75,8.25) [circ] {};
					\node at (10,8.25) [circ] {};
					\node at (11.25,8.25) [circ] {};
					\node at (12.5,8.25) [circ] {};
					\node at (13.75,8.25) [circ] {};
					\node at (-3.75,9.5) [circ] {};
					\node at (-2.5,9.5) [circ] {};
					\node at (-1.25,9.5) [circ] {};
					\node at (0,9.5) [circ] {};
					\node at (1.25,9.5) [circ] {};
					\node at (2.5,9.5) [circ] {};
					\node at (3.75,9.5) [circ] {};
					\node at (5,9.5) [circ] {};
					\node at (6.25,9.5) [circ] {};
					\node at (7.5,9.5) [circ] {};
					\node at (8.75,9.5) [circ] {};
					\node at (10,9.5) [circ] {};
					\node at (11.25,9.5) [circ] {};
					\node at (12.5,9.5) [circ] {};
					\node at (-1.25,10.75) [circ] {};
					\node at (0,10.75) [circ] {};
					\node at (1.25,10.75) [circ] {};
					\node at (2.5,10.75) [circ] {};
					\node at (3.75,10.75) [circ] {};
					\node at (5,10.75) [circ] {};
					\node at (6.25,10.75) [circ] {};
					\node at (7.5,10.75) [circ] {};
					\node at (8.75,10.75) [circ] {};
					\node at (10,10.75) [circ] {};
					\node at (-1.25,12) [circ] {};
					\node at (0,12) [circ] {};
					\node at (1.25,12) [circ] {};
					\node at (2.5,12) [circ] {};
					\node at (3.75,12) [circ] {};
					\node at (5,12) [circ] {};
					\node at (6.25,12) [circ] {};
					\node at (7.5,12) [circ] {};
					\node at (8.75,12) [circ] {};
					\node at (10,12) [circ] {};
					\node at (11.25,12) [circ] {};
					\node at (-3.75,13.25) [circ] {};
					\node at (-2.5,13.25) [circ] {};
					\node at (-1.25,13.25) [circ] {};
					\node at (0,13.25) [circ] {};
					\node at (1.25,13.25) [circ] {};
					\node at (2.5,13.25) [circ] {};
					\node at (3.75,13.25) [circ] {};
					\node at (5,13.25) [circ] {};
					\node at (6.25,13.25) [circ] {};
					\node at (7.5,13.25) [circ] {};
					\node at (8.75,13.25) [circ] {};
					\node at (10,13.25) [circ] {};
					\node at (-1.25,14.5) [circ] {};
					\node at (0,14.5) [circ] {};
					\node at (1.25,14.5) [circ] {};
					\node at (2.5,14.5) [circ] {};
					\node at (3.75,14.5) [circ] {};
					\node at (5,14.5) [circ] {};
					\node at (6.25,14.5) [circ] {};
					\node at (7.5,14.5) [circ] {};
					\node at (8.75,14.5) [circ] {};
					\node at (10,14.5) [circ] {};
					\node at (11.25,14.5) [circ] {};
					\node at (12.5,14.5) [circ] {};
					\node at (2.5,15.75) [circ] {};
					\node at (3.75,15.75) [circ] {};
					\node at (5,15.75) [circ] {};
					\node at (6.25,15.75) [circ] {};
					\node at (7.5,15.75) [circ] {};
					\node at (10,15.75) [circ] {};
					\node at (11.25,15.75) [circ] {};
					\draw [line width=1pt, short] (7,4.25) -- (7,4.25);
					\draw [line width=1pt, short] (7.5,4.5) -- (8.75,3.25);
					\node at (2.5,17) [circ] {};
					\node at (3.75,17) [circ] {};
					\node at (5,17) [circ] {};
					\node at (6.25,17) [circ] {};
					\node at (7.5,17) [circ] {};
					\draw [line width=1pt, short] (0,5.75) .. controls (1.25,6.25) and (1.25,6.5) .. (2.5,5.75);
					\draw [line width=1pt, short] (1.25,5.75) -- (2.5,5.75);
					\draw [line width=1pt, short] (2.5,5.75) -- (1.25,4.5);
					\draw [line width=1pt, short] (2.5,5.75) -- (3.75,5.75);
					\draw [line width=1pt, short] (2.5,4.5) -- (3.75,5.75);
					\draw [line width=1pt, short] (3.75,5.75) -- (5,4.5);
					\draw [line width=1pt, short] (3.75,5.75) -- (5,5.75);
					\draw [line width=1pt, short] (5,5.75) -- (6.25,4.5);
					\draw [line width=1pt, short] (5,5.75) -- (6.25,5.75);
					\draw [line width=1pt, short] (6.25,5.75) -- (7.5,4.5);
					\draw [line width=1pt, short] (6.25,5.75) -- (7.5,5.75);
					\draw [line width=1pt, short] (7.5,5.75) -- (8.75,5.75);
					\draw [line width=1pt, short] (7.5,4.5) -- (8.75,4.5);
					\draw [line width=1pt, short] (10,5.75) -- (15,5.75);
					\draw [line width=1pt, short] (10,5.75) -- (11.25,4.5);
					\draw [line width=1pt, short] (10,5.75) -- (12.5,4.5);
					\node at (15,4.5) [circ] {};
					\draw [line width=1pt, short] (12.5,5.75) -- (13.75,4.5);
					\draw [line width=1pt, short] (12.5,5.75) -- (15,4.5);
					\draw [line width=1pt, short] (3.75,5.75) -- (3.75,7);
					\draw [line width=1pt, short] (3.75,7) -- (3.75,17);
					\draw [line width=1pt, short] (3.75,7) -- (5,7);
					\draw [line width=1pt, short] (3.75,8.25) -- (5,8.25);
					\draw [line width=1pt, short] (3.75,9.5) -- (5,9.5);
					\draw [line width=1pt, short] (3.75,9.5) -- (5,10.75);
					\draw [line width=1pt, short] (3.75,12) -- (5,12);
					\draw [line width=1pt, short] (3.75,13.25) -- (5,13.25);
					\draw [line width=1pt, short] (2.5,17) -- (5,17);
					\draw [line width=1pt, short] (2.5,15.75) -- (3.75,15.75);
					\draw [line width=1pt, short] (2.5,14.5) -- (3.75,14.5);
					\draw [line width=1pt, short] (3.75,12) -- (2.5,12);
					\draw [line width=1pt, short] (2.5,13.25) -- (2.5,12);
					\draw [line width=1pt, short] (1.25,12) -- (2.5,12);
					\draw [line width=1pt, short] (1.25,13.25) -- (2.5,12);
					\draw [line width=1pt, short] (2.5,10.75) -- (3.75,10.75);
					\draw [line width=1pt, short] (2.5,8.25) -- (3.75,8.25);
					\draw [line width=1pt, short] (2.5,8.25) -- (2.5,7);
					\draw [line width=1pt, short] (2.5,8.25) -- (1.25,7);
					\draw [line width=1pt, short] (2.5,8.25) -- (-3.75,8.25);
					\draw [line width=1pt, short] (0,7) -- (1.25,7);
					\draw [line width=1pt, short] (0,8.25) -- (-1.25,7);
					\draw [line width=1pt, short] (-1.25,8.25) -- (-2.5,7);
					\draw [line width=1pt, short] (-2.5,8.25) -- (-3.75,7);
					\draw [line width=1pt, short] (-3.75,9.5) -- (-1.25,8.25);
					\draw [line width=1pt, short] (-2.5,9.5) -- (-1.25,8.25);
					\draw [line width=1pt, short] (-1.25,9.5) -- (0,8.25);
					\draw [line width=1pt, short] (1.25,9.5) -- (2.5,8.25);
					\draw [line width=1pt, short] (0,9.5) -- (2.5,9.5);
					\draw [line width=1pt, short] (0,12) -- (0,10.75);
					\draw [line width=1pt, short] (-1.25,10.75) -- (0,10.75);
					\draw [line width=1pt, short] (-1.25,12) -- (0,10.75);
					\draw [line width=1pt, short] (0,10.75) -- (1.25,10.75);
					\draw [line width=1pt, short] (-2.5,13.25) -- (0,13.25);
					\draw [line width=1pt, short] (0,14.5) -- (0,13.25);
					\draw [line width=1pt, short] (-1.25,14.5) -- (0,13.25);
					\draw [line width=1pt, short] (1.25,14.5) -- (0,13.25);
					\draw [line width=1pt, short] (5,5.75) -- (6.25,7);
					\draw [line width=1pt, short] (6.25,5.75) -- (7.5,7);
					\draw [line width=1pt, short] (7.5,5.75) -- (8.75,7);
					\draw [line width=1pt, short] (7.5,5.75) -- (10,7);
					\draw [line width=1pt, short] (12.5,5.75) -- (13.75,7);
					\draw [line width=1pt, short] (13.75,5.75) -- (15,7);
					\draw [line width=1pt, short] (8.75,7) -- (10,8.25);
					\draw [line width=1pt, short] (8.75,8.25) -- (8.75,7);
					\draw [line width=1pt, short] (10,8.25) -- (12.5,8.25);
					\draw [line width=1pt, short] (10,9.5) -- (10,8.25);
					\draw [line width=1pt, short] (10,8.25) -- (11.25,9.5);
					\draw [line width=1pt, short] (11.25,8.25) -- (12.5,9.5);
					\draw [line width=1pt, short] (10,8.25) -- (11.25,7);
					\draw [line width=1pt, short] (11.25,8.25) -- (12.5,7);
					\draw [line width=1pt, short] (6.25,8.25) -- (7.5,8.25);
					\draw [line width=1pt, short] (6.25,8.25) -- (8.75,10.75);
					\draw [line width=1pt, short] (7.5,9.5) -- (7.5,10.75);
					\draw [line width=1pt, short] (6.25,10.75) -- (7.5,10.75);
					\draw [line width=1pt, short] (6.25,9.5) -- (7.5,10.75);
					\draw [line width=1pt, short] (6.25,12) -- (7.5,10.75);
					\draw [line width=1pt, short] (8.75,12) -- (8.75,10.75);
					\draw [line width=1pt, short] (7.5,12) -- (10,12);
					\draw [line width=1pt, short] (8.75,12) -- (10,10.75);
					\draw [line width=1pt, short] (6.25,13.25) -- (6.25,17);
					\draw [line width=1pt, short] (6.25,15.75) -- (7.5,17);
					\draw [line width=1pt, short] (6.25,14.5) -- (7.5,15.75);
					\draw [line width=1pt, short] (5,15.75) -- (6.25,14.5);
					\draw [line width=1pt, short] (5,14.5) -- (6.25,13.25);
					\draw [line width=1pt, short] (7.5,14.5) -- (6.25,13.25);
					\draw [line width=1pt, short] (6.25,13.25) -- (7.5,13.25);
					\draw [line width=1pt, short] (7.5,13.25) -- (8.75,14.5);
					\draw [line width=1pt, short] (7.5,13.25) -- (10,14.5);
					\draw [line width=1pt, short] (7.5,13.25) -- (8.75,13.25);
					\draw [line width=1pt, short] (10,15.75) -- (10,14.5);
					\draw [line width=1pt, short] (10,14.5) -- (11.25,14.5);
					\draw [line width=1pt, short] (10,14.5) -- (11.25,15.75);
					\draw [line width=1pt, short] (-2.5,8.25) .. controls (-3.75,8.75) and (-3.5,8.75) .. (-5,8.25);
					\draw [line width=1pt, short] (-3.75,13.25) .. controls (-2.25,14) and (-2.5,14) .. (-1.25,13.25);
					\draw [line width=1pt, short] (0,13.25) .. controls (2.25,14.25) and (2.25,14.25) .. (3.75,13.25);
					\draw [line width=1pt, short] (3.75,13.25) .. controls (5.25,14) and (5,13.75) .. (6.25,13.25);
					\draw [line width=1pt, short] (7.5,13.25) .. controls (9,12.75) and (9,12.5) .. (10,13.25);
					\draw [line width=1pt, short] (10,14.5) .. controls (11.5,15.25) and (11.5,15) .. (12.5,14.5);
					\draw [line width=1pt, short] (8.75,12) .. controls (10,11.25) and (10.25,11.5) .. (11.25,12);
					\draw [line width=1pt, short] (8.75,9.5) .. controls (8.25,8.25) and (8,8.25) .. (8.75,7);
					\draw [line width=1pt, short] (11.25,8.25) .. controls (12.75,9) and (12.5,9) .. (13.75,8.25);
					\draw [line width=1pt, short] (13.75,5.75) .. controls (15.25,6.5) and (15.25,6.25) .. (16.25,5.75);
					\draw [line width=1pt, short] (7.5,5.75) .. controls (8.75,5.25) and (9,5.25) .. (10,5.75);
					\draw [line width=1pt, short] (7.5,4.5) .. controls (9,3.75) and (9,4) .. (10,4.5);
					\draw [line width=1pt, short] (3.75,8.25) .. controls (5.25,9) and (5.25,8.75) .. (6.25,8.25);
					\draw [line width=1pt, short] (1.25,10.75) .. controls (2.5,10.25) and (2.5,10) .. (3.75,10.75);
					\node [font=\Large] at (3.75,5.25) {2};
					\node [font=\Large] at (5,4.25) {};
					\node [font=\Large] at (5,4.2) {3};
					\node [font=\Large] at (2.5,4.25) {};
					\node [font=\Large] at (2.5,4.2) {6};
					\node [font=\Large] at (1.25,4.15) {67};
					\node [font=\Large] at (6.25,4.25) {};
					\node [font=\Large] at (6.25,4.2) {8};
					\node [font=\Large] at (7.35,4.2) {78};
					\node [font=\Large] at (9.25,3.25) {79};
					\node [font=\Large] at (9,4.75) {76};
					\node [font=\Large] at (10,4.25) {};
					\node [font=\Large] at (10,4) {74};
					\node [font=\Large] at (11.25,4.25) {60};
					\node [font=\Large] at (12.5,4.25) {63};
					\node [font=\Large] at (0,5.5) {70};
					\node [font=\Large] at (1.25,5.5) {64};
					\node [font=\Large] at (2.5,5.5) {};
					\node [font=\Large] at (2.5,5.25) {66};
					\node [font=\Large] at (5,5.25) {10};
					\node [font=\Large] at (6,6) {14};
					\node [font=\Large] at (7.5,5.3) {46};
					\node [font=\Large] at (8.75,6) {38};
					\node [font=\Large] at (10,5.25) {62};
					\node [font=\Large] at (12,6) {126};
					\node [font=\Large] at (0,6.75) {120};
					\node [font=\Large] at (1.25,6.75) {112};
					\node [font=\Large] at (2.5,6.75) {56};
					\node [font=\Large] at (5,6.75) {98};
					\node [font=\Large] at (4,7.25) {34};
					\node [font=\Large] at (3.5,8) {32};
					\node [font=\Large] at (5,8) {40};
					\node [font=\Large] at (6.25,8) {96};
					\node [font=\Large] at (7.5,8) {97};
					\node [font=\Large] at (6.2,7.25) {42};
					\node [font=\Large] at (7.45,7.25) {12};
					\node [font=\Large] at (8.5,7) {};
					\node [font=\Large] at (8.25,7) {  47};
					\node [font=\Large] at (10,6.75) {44};
					\node [font=\Large] at (11,6.75) {7};
					\node [font=\Large] at (12,7) {111};
					\node [font=\Large] at (12.5,8) {};
					\node [font=\Large] at (10.25,9.75) {55};
					\node [font=\Large] at (11.5,9.75) {35};
					\node [font=\Large] at (12.5,9.75) {102};
					\node [font=\Large] at (9.1,8.25) {15};
					\node [font=\Large] at (9,9.25) {45};
					\node [font=\Large] at (7.5,9) {104};
					\node [font=\Large] at (6.25,9.25) {124};
					\node [font=\Large] at (5,9.25) {41};
					\node [font=\Large] at (5,10.25) {49};
					\node [font=\Large] at (6,10.5) {109};
					\node [font=\Large] at (7.95,10.9) {108};
					\node [font=\Large] at (9.1,10.5) {105};
					\node [font=\Large] at (10.5,10.75) {106};
					\node [font=\Large] at (11.25,11.55) {123};
					\node [font=\Large] at (10.25,12.25) {43};
					\node [font=\Large] at (8.75,12.25) {107};
					\node [font=\Large] at (7.5,12.25) {99};
					\node [font=\Large] at (6,12.25) {100};
					\node [font=\Large] at (5,11.75) {11};
					\node [font=\Large] at (4,11.75) {9};
					\node [font=\Large] at (2.5,11.70) {73};
					\node [font=\Large] at (1.25,11.70) {75};
					\node [font=\Large] at (1,13) {72};
					\node [font=\Large] at (2,13.25) {  89};
					\node [font=\Large] at (3.25,13) {   25};
					\node [font=\Large] at (0,12.25) {65};
					\node [font=\Large] at (-1.25,12.25) {77};
					\node [font=\Large] at (-1.5,11) {101};
					\node [font=\Large] at (0,10.5) {69};
					\node [font=\Large] at (1.25,11) {5};
					\node [font=\Large] at (2.25,11) {0};
					\node [font=\Large] at (3.5,11) {1};
					\node [font=\Large] at (0.25,9.25) {51};
					\node [font=\Large] at (1.25,9.75) {50};
					\node [font=\Large] at (2.25,9.75) {114};
					\node [font=\Large] at (3.35,9.5) {33};
					\node [font=\Large] at (1.8,8.5) {48};
					\node [font=\Large] at (0.25,8) {52};
					\node [font=\Large] at (-1,8) {20};
					\node [font=\Large] at (-1,9.75) {53};
					\node [font=\Large] at (-2.25,9.75) {4};
					\node [font=\Large] at (4.75,13) {24};
					\node [font=\Large] at (6,13) {27};
					\node [font=\Large] at (7,13) { 19};
					\node [font=\Large] at (9,13.5) {17};
					\node [font=\Large] at (10.35,13.35) {18};
					\node [font=\Large] at (10,14.15) {83};
					\node [font=\Large] at (11.25,14) {87};
					\node [font=\Large] at (12.5,14.25) {115};
					\node [font=\Large] at (11.75,15.75) {82};
					\node [font=\Large] at (10.25,16) {81};
					\node [font=\Large] at (9,14.75) {23};
					\node [font=\Large] at (7.75,14.75) {59};
					\node [font=\Large] at (7.75,16) {30};
					\node [font=\Large] at (7.95,17) {88};
					\node [font=\Large] at (6.5,17.25) {91};
					\node [font=\Large] at (5,17.25) {117};
					\node [font=\Large] at (3.75,17.25) {125};
					\node [font=\Large] at (2.25,17.25) {127};
					\node [font=\Large] at (2,16) {113};
					\node [font=\Large] at (3.35,16.1) {121};
					\node [font=\Large] at (2.25,14.75) {61};
					\node [font=\Large] at (3.5,14.75) {57};
					\node [font=\Large] at (4.75,14.5) {31};
					\node [font=\Large] at (4.5,15.75) {  58};
					\node [font=\Large] at (1.25,14.75) {28};
					\node [font=\Large] at (0,14.75) {21};
					\node [font=\Large] at (-1.25,14.75) {13};
					\node [font=\Large] at (0,13) {29};
					\node [font=\Large] at (-1.25,13) {93};
					\node [font=\Large] at (-2.5,13) {95};
					\node [font=\Large] at (-4,13) {92};
					\node [font=\Large] at (-4,9.75) {16};
					\node [font=\Large] at (-5,8) {85};
					\node [font=\Large] at (-3.75,8) {80};
					\node [font=\Large] at (-3.75,6.75) {68};
					\node [font=\Large] at (-2.5,6.75) {22};
					\node [font=\Large] at (-1.25,6.75) {36};
					\node [font=\Large] at (-2.3,7.9) {84};
					\node [font=\Large] at (13.75,4.25) {110};
					\node [font=\Large] at (15,4.25) {94};
					\node [font=\Large] at (14,5.5) {118};
					\node [font=\Large] at (15.25,5.5) {86};
					\node [font=\Large] at (16.75,5.75) {116};
					\node [font=\Large] at (13.75,7.25) {122};
					\node [font=\Large] at (15,7.25) {54};
					\node [font=\Large] at (12.5,8) {};
					\node [font=\Large] at (12.5,7.75) {71};
					\node [font=\Large] at (14.25,8.5) {119};
					\node [font=\Large] at (10,7.75) {39};
					\node [font=\Large] at (11,8.55) {103};
					\node [font=\Large] at (5.85,16) {90};
					\node [font=\Large] at (5.75,14.5) {26};
					\node at (-2.5,10.75) [circ] {};
					\draw [line width=1pt, short] (1.25,10.75) .. controls (0,9.75) and (-0.5,10.25) .. (-2.5,10.75);
					\node [font=\Large] at (-3,10.5) {37};
				\end{circuitikz}
			}%
			
				\caption {$T_3$: Third CIST in Q7 with diameter 17}
		\end{figure}
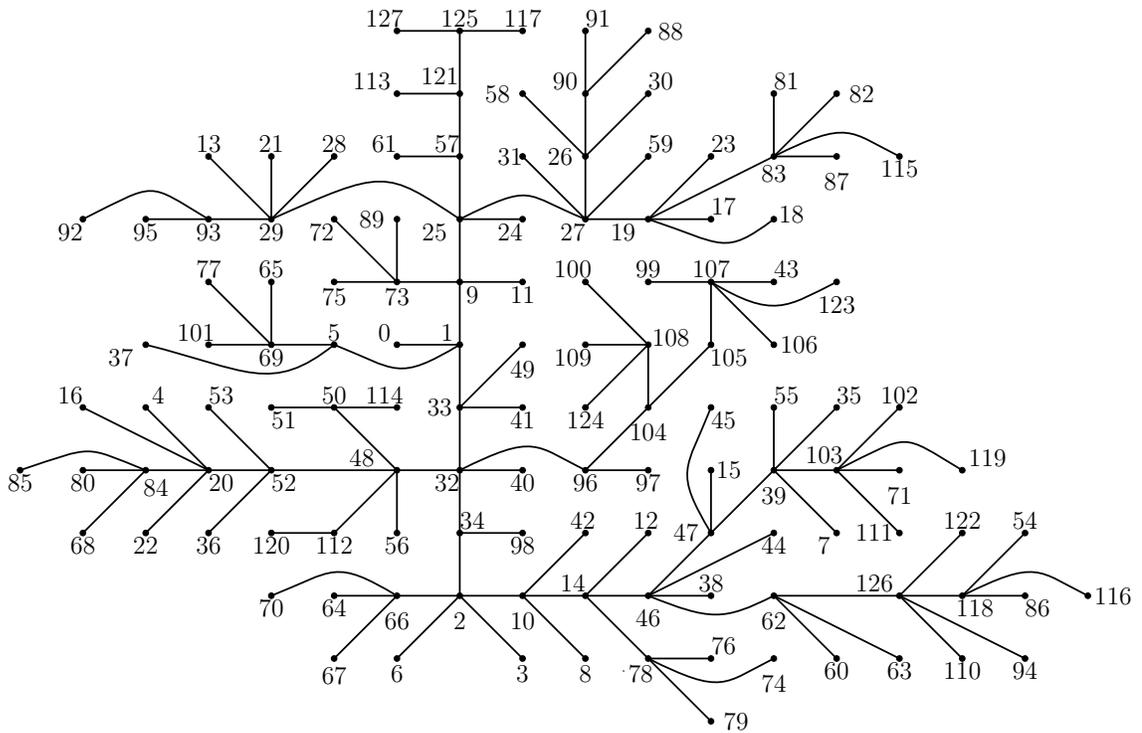
	

 	Here, we observe that\\
		
		\noindent
		$InV(T_1)= \{0 , 6 , 7 , 8 , 15 , 16 , 17 , 21 , 23 , 24 , 31 , 35 , 38 , 40 , 44 , 45 , 51 , 54 , 55 , 59 , 63 , 64 , 71 , 72 , 74 ,\\ 77 , 79 , 80,  82 , 85, 88, 91, 92 , 97 , 99 , 100 , 102 , 110 , 111 , 113 , 115 , 122 , 123 , 124\}$,\\
		
		\noindent
		$InV(T_2) = \{3, 4 , 11 , 12 , 13 , 18 ,22 , 28, 30 , 36 , 37 , 41 , 42 , 43 , 49 , 53 , 56 , 58 , 60 , 61 , 65 , 67 , 68 , 70, 75 ,\\ 76 , 81 , 86 , 87 , 89 , 94 , 95 , 98 , 101 , 106 , 109 , 114 , 117 , 119 , 120 , 127 \}$ \\
		
		\noindent
		$InV(T_3)= \{1 , 2 , 5 , 9 , 10 , 14 , 19 , 20 , 25 , 26 , 27 , 29 , 32 , 33 , 34 , 39 , 46 , 47 , 48 , 50 , 52 , 57 , 62 , 66 , 69 ,\\73 , 78 , 83 , 84 , 90 , 93 , 96 , 103 , 104 , 105 , 107 , 108 , 112 , 118 , 121 , 125 , 126\}$	\\

		are such that $InV(T_i) \cap InV(T_j) = \phi $, for $i \neq j$ and $ 1 \leq i, j \leq 3.$ Also, note that the edges of the trees above are listed below: \\
		
		Edge set of first CIST $T_1:$ \\
		
		$\{ \langle0 , 2  \rangle,   \langle 0, 8 \rangle,    \langle0 , 16 \rangle,    \langle0 , 32 \rangle,   \langle0 , 64 \rangle ,   \langle1 , 17\rangle,   \langle3 ,35 \rangle,    \langle4 , 6 \rangle,  \langle5 , 21 \rangle ,   \langle6 , 7\rangle,    \langle6 , 22  \rangle,   \langle7 , 15 \rangle,     \langle7 , 71 \rangle, \\   \langle8 , 9 \rangle,    \langle8 , 40 \rangle,     \langle8 , 72  \rangle,     \langle10 , 74 \rangle,    \langle11 , 15 \rangle,    \langle12 , 44\rangle,    \langle13 , 45\rangle,    \langle14 , 15\rangle,    \langle15 , 31 \rangle,    \langle15 , 79 \rangle,     \langle16 , 17 \rangle, \\   \langle16 , 24 \rangle,    \langle16 , 48 \rangle,   \langle16 , 80 \rangle,        \langle17 , 21 \rangle,    \langle17 , 25 \rangle,  \langle18 , 82 \rangle,  \langle19 , 51 \rangle,       \langle20 , 21\rangle,    \langle21 , 23 \rangle,    \langle21 , 85 \rangle,  \\  \langle 23, 31 \rangle,   \langle23 , 55 \rangle,    \langle24 , 26 \rangle,    \langle24 , 28 \rangle,  \langle24 , 88 \rangle,    \langle27 , 91\rangle,   \langle29 , 31 \rangle,    \langle30 , 31 \rangle,    \langle33 , 35 \rangle,     \langle34 , 38 \rangle,    \langle35 , 51\rangle,  \\  \langle35 , 99 \rangle,     \langle36 , 44 \rangle,    \langle37 , 45 \rangle,     \langle38 , 39 \rangle,    \langle38 , 54 \rangle,   \langle38 , 102  \rangle,    \langle40 , 42 \rangle,    \langle40 , 44 \rangle,    \langle40 , 56\rangle,     \langle41 , 45 \rangle, \\   \langle43 , 59 \rangle,    \langle44 , 45\rangle,    \langle45 , 61\rangle,  \langle 46 , 110 \rangle,  \langle47 , 63\rangle,     \langle49 , 113 \rangle,    \langle50 , 54 \rangle,     \langle51 , 55 \rangle,   \langle51 , 59 \rangle,     \langle51 , 115\rangle,    \langle52 , 54 \rangle, \\   \langle53 , 55\rangle,    \langle54 , 55\rangle,     \langle54 , 62\rangle,    \langle57 , 59\rangle,    \langle58 , 122\rangle,    \langle59 , 63\rangle,    \langle59 , 123\rangle,    \langle60 , 124\rangle,   \langle63 , 127\rangle,    \langle64 , 65\rangle,     \langle66 , 82\rangle, \\    \langle67 , 71\rangle,   \langle68 , 100\rangle,    \langle69 , 71\rangle, \langle 70 , 71\rangle,   \langle72 , 74\rangle,     \langle72 , 104 \rangle,   \langle73 , 77\rangle,    \langle75 , 79\rangle,    \langle76 , 77\rangle,     \langle77 , 79\rangle,    \langle78 , 110\rangle, \\   \langle79 , 111\rangle,    \langle80 , 82\rangle,     \langle81 , 85\rangle,     \langle82 , 86\rangle,     \langle82 , 90\rangle,     \langle83 , 91\rangle,    \langle84 , 92\rangle,     \langle85 , 87\rangle,     \langle85 , 93 \rangle,    \langle88 , 92\rangle,    \langle88 , 120\rangle, \\   \langle89 , 91\rangle,    \langle91 , 95\rangle,     \langle91 ,     123 \rangle,    \langle92 , 94\rangle,     \langle92 , 
		124 \rangle,   \langle96 , 100\rangle,    \langle97 , 99\rangle,     \langle97 , 105 \rangle,    \langle98 , 99\rangle,     \langle99 , 103\rangle,     \langle100 , 101\rangle, \\  \langle100 , 102 \rangle,    \langle100 , 116 \rangle,  \langle 102 , 118 \rangle,  \langle106 , 122\rangle,    \langle107 , 111\rangle,    \langle108 , 110\rangle,    \langle109 , 111\rangle,    \langle110 , 111\rangle,    \langle112 , 113\rangle,  \\   \langle113 , 115\rangle,   \langle113 , 
		117 \rangle,     \langle114 , 122\rangle,   \langle115 , 119\rangle,     \langle121 , 123\rangle,   \langle122 , 123\rangle,     \langle124 , 125\rangle,    \langle124 , 126\rangle \}$\\
		
		
		Edge set of second CIST $T_2:$\\
		
		$\{ \langle0 , 4\rangle,    \langle1 , 3 \rangle,
		\langle2 , 18\rangle,    \langle3 , 7\rangle, 
		\langle3 , 11 \rangle,    \langle3 ,  19\rangle,
		\langle3 , 67 \rangle,   \langle4 , 5 \rangle,       \langle4 ,  12\rangle,    \langle4 , 68 \rangle, 		
		\langle 6 , 70 \rangle,  
		\langle8 ,  12\rangle, \\   \langle9 , 13 \rangle,  
		\langle10 , 11\rangle,    \langle11 , 27 \rangle,   \langle11 , 43 \rangle,   \langle11 , 75\rangle,
		\langle12 ,  13\rangle,    \langle12 , 28\rangle,
		\langle12 , 76 \rangle,     \langle13 , 15\rangle,
		\langle14 , 30 \rangle,    \langle16 , 18 \rangle, \\
		\langle17 , 49\rangle,    \langle18 , 22 \rangle, 
		\langle18 , 26 \rangle,    \langle20 , 28 \rangle,
		\langle21 , 53 \rangle,    \langle22 , 23\rangle,   \langle22 , 30 \rangle,    \langle22 , 54 \rangle, 
		\langle22 , 86\rangle,    \langle24 , 56 \rangle,
		\langle25 , 89\rangle, \\   \langle28 , 30 \rangle, 
		\langle28 , 60\rangle,     \langle29 , 61 \rangle,
		\langle30 , 94 \rangle,     \langle31 , 95 \rangle, 
		\langle32 , 36 \rangle,    \langle33 , 37 \rangle,
		\langle34 , 42 \rangle,    \langle35 , 43 \rangle,
		\langle36 , 37\rangle,     \langle36 , 38 \rangle, \\ \langle36 ,  100\rangle,     \langle37 ,  39\rangle, \langle37 , 53\rangle,      \langle37 , 101 \rangle,
		\langle40 , 41 \rangle,     \langle41 , 43 \rangle, \langle41 ,  105\rangle,     \langle42 , 43 \rangle, \langle42 , 46 \rangle,      \langle42 ,  58\rangle, \\ \langle42 , 106 \rangle,     \langle43 , 47\rangle,
		\langle44 , 60 \rangle,     \langle45 , 109 \rangle,
		\langle48 , 49 \rangle,     \langle49 , 51 \rangle,
		\langle49 , 53 \rangle,     \langle49 , 57 \rangle,    \langle50 , 58\rangle,  
		\langle52 , 60 \rangle,  \\   \langle53 , 61 \rangle, \langle53 , 117 \rangle,    \langle55 , 119 \rangle,
		    \langle56 , 58 \rangle, \langle56 , 60 \rangle,  \langle56 , 120 \rangle,     \langle58 , 59 \rangle, \langle58 , 62 \rangle,      \langle60 , 61 \rangle, 
		\langle61 , 63\rangle,  \\    \langle64 , 68\rangle,
		\langle65 , 67\rangle,      \langle65 , 73 \rangle, \langle65 , 81 \rangle,     \langle66 ,  98\rangle,
		\langle67 , 83 \rangle,     \langle67 , 99 \rangle,
		\langle68 , 69 \rangle,     \langle68 ,  70\rangle,
		\langle71 , 87 \rangle, \\    \langle72 , 76 \rangle,
		\langle74 ,  75\rangle,   \langle75 , 91 \rangle,   \langle75 , 107 \rangle,    \langle76 , 92 \rangle, 
		\langle76 , 108 \rangle,    \langle77 , 109 \rangle,
		\langle78 , 94 \rangle,     \langle79 , 95 \rangle, 
		\langle80 ,  81\rangle,     \langle81 , 89 \rangle, \\	\langle 81 , 113 \rangle,   \langle82 , 114 \rangle,    \langle84 , 86 \rangle,
		\langle85 , 117 \rangle,    \langle86 , 87 \rangle,
		\langle88 , 89 \rangle,     \langle89 , 93 \rangle,
		\langle90 , 94\rangle,     \langle94 , 95 \rangle,
		\langle95 , 127 \rangle,    \langle96 , 98 \rangle, \\
		\langle97 , 101 \rangle,     \langle98 , 102\rangle, \langle98 , 106 \rangle,    \langle98 , 114 \rangle,
		\langle101 , 103\rangle,  \langle101 , 109 \rangle,  \langle104 , 106 \rangle,
		\langle106 , 110 \rangle,    
		\langle109, 125\rangle, \\    \langle111 , 127\rangle,
		\langle112 , 114 \rangle,     \langle114 , 115\rangle,
		\langle116 , 117 \rangle,    \langle118 , 119 \rangle,
		\langle119 , 127 \rangle,    \langle120 , 121 \rangle,  \langle120 , 124 \rangle, \\    \langle120 , 122 \rangle,
		\langle123 , 127 \rangle,    \langle126 , 127\rangle \}$ \\

		
		Edge set of third CIST $T_3:$\\
		
		$\{ \langle0 , 1\rangle,     \langle1 , 5\rangle,
		\langle 1 , 9\rangle,       \langle1 , 33\rangle,
		\langle2 , 3 \rangle,       \langle2 , 6\rangle, 
		\langle 2 , 10\rangle,      \langle2 , 34 \rangle, 
		\langle2 , 66\rangle,       \langle4 , 20\rangle,
		\langle5 , 37\rangle,      \langle5 , 69\rangle,
		\langle7 , 39\rangle,  \\    \langle8 , 10\rangle,
		\langle9 , 11\rangle,      \langle9 , 25\rangle, 
		\langle9 , 73\rangle,     \langle10 , 14\rangle, 
		\langle10 , 42\rangle,    \langle12 , 14\rangle,
		\langle13 , 29\rangle,    \langle14 , 46\rangle, 
		\langle14 , 78\rangle,    \langle15 , 47\rangle, \\
		\langle16 , 20\rangle,    \langle17 , 19\rangle,
		\langle18 , 19\rangle,    \langle19 , 23\rangle, 
		\langle19 , 27\rangle,    \langle19 , 83\rangle,
		\langle20 , 22\rangle,   
		\langle 20 , 52\rangle,    \langle 20, 84\rangle,
		\langle21 , 29\rangle,    \langle24 , 25\rangle,\\
		\langle 25 , 27\rangle,    \langle25 , 29\rangle, 
		\langle 25 , 57\rangle,    \langle26 , 27\rangle, 
		\langle26 , 30\rangle,     \langle26, 58\rangle, 
		\langle26 , 90\rangle,     \langle27 , 59\rangle, 
		\langle27 , 31\rangle,      \langle 28 , 29\rangle, 
		\langle29 , 93\rangle, \\
		\langle32 , 33\rangle,     \langle32 , 34 \rangle, 
		\langle32 , 40\rangle,     \langle32 , 48\rangle, 
		\langle32 , 96\rangle,     \langle33, 41\rangle, 
		\langle33 , 49\rangle,     \langle34 , 98\rangle,
		\langle35 , 39\rangle,    \langle36 , 52\rangle,
		\langle38 , 46\rangle, \\   \langle39 , 47\rangle, 
		\langle39 , 55\rangle,    \langle39 , 103\rangle,
		\langle43 , 107\rangle,   \langle44 , 46\rangle,
		\langle45 , 47\rangle,   \langle46 , 47\rangle, 
		\langle46 , 62\rangle,    \langle48 , 50\rangle, 
		\langle48 , 52\rangle, \\   \langle48 , 56\rangle, 
		\langle 48 , 112\rangle,   \langle50 , 51\rangle, 
		\langle50 , 114\rangle,    \langle52 , 53\rangle,
		\langle54 , 118\rangle,    \langle57 , 61\rangle, 
		\langle 57 , 121\rangle,    \langle60 , 62\rangle,
		\langle62 , 63\rangle,  \\    \langle62 , 126\rangle,
		\langle64 , 66\rangle,     \langle65 , 69\rangle,
		\langle66 , 67\rangle,     \langle66 , 70\rangle,
		\langle68 , 84\rangle,     \langle69 , 77\rangle, 
		\langle69 , 101\rangle,    \langle71 , 103\rangle,
		\langle72 , 73\rangle,  \\   \langle73 , 75\rangle, 
		\langle 73 , 89\rangle,    \langle74 , 78\rangle,
		\langle76 , 78\rangle,     \langle 78 , 79\rangle,
		\langle80 , 84\rangle,     \langle81 , 83\rangle,
		\langle82 , 83\rangle,     \langle83 , 87\rangle,
		\langle83 , 115\rangle,    \langle84 , 85\rangle, \\
		\langle86 , 118\rangle,    \langle88 , 90\rangle,
		\langle90 , 91\rangle,     \langle92 , 93\rangle,
		\langle93 , 95\rangle,     \langle94 , 126\rangle,
		\langle96 , 97\rangle,      \langle96 , 104\rangle,    
		\langle99 , 107\rangle,    \langle100, 108\rangle, \\
		\langle102 , 103\rangle,   \langle103 , 111\rangle, 
		\langle 103, 119\rangle,   \langle104 , 105\rangle, 
	    \langle104 , 108\rangle,     \langle105 , 107\rangle,
		\langle106 , 107\rangle,    \langle107 , 123\rangle, \\
		\langle108 , 109\rangle,    \langle108 , 124\rangle,
		\langle110 , 126\rangle,    \langle112 , 120\rangle,
		\langle113 , 121\rangle,    \langle116 , 118\rangle,
		\langle117 , 125\rangle,    \langle118 , 126\rangle, \\
		\langle121 , 125\rangle,     \langle122, 126\rangle,
		\langle125 , 127\rangle \}$\\
		
			are such that $E(T_i) \cap E(T_j) = \phi $, for $i \neq j$ and $ 1 \leq i, j \leq 3.$\\
			
			We will now refer to the following theorem\cite{h1}, which allows us to verify that the trees \( T_1 \), \( T_2 \), and \( T_3 \) constructed above are CISTs on \( Q_7 \).   
		
		\begin{thm} [\cite{h1} \label{lm1}] Let $k \geq 2$ be an integer. $T_1, T_2,\ldots ,T_k$ are completely independent spanning trees in a graph $G$ if and only if they are edge-disjoint spanning trees of $G$ and for any $v \in V(G)$, there is at most one $T_i$ such that $d_{T_i}(v) > 1$.
		\end{thm}

		\section {\textbf{Construction of  three CISTs in hypercubes} }
		In this section, we construct three CISTs in the hypercube \( Q_n \), where \( n \geq 8 \),  using induction. Specifically, among these spanning trees, two of these trees have diameters of at most \(2n+1\) and \(2n+3\), respectively, and the third tree has a diameter of at most \(2n+4\).
		
		To proceed, we first introduce the necessary definitions and results.

We shall employ the following corollary \cite{sm}.
		
		\begin{cor} [\cite{sm}] \label{lm3}
			Let $G_n$ be the $n-$dimensional variant hypercube for $n \geq 4$ and suppose that $T_i$, for $1 \leq i \leq k $ $(k < n)$ be $k$ CISTs of $G_n$. Let $\overline T_i$ be a spanning tree of $G_{n+1}$ constructed from $T^0_i$ and $T^1_i$ by adding an edge $\langle u_i, v_i \rangle \in E(G_{n+1})$ to connect two internal vertices $u_i \in V(T^0_i)$ and $v_i \in V(T^1_i)$. Then, $\overline T_i$ are $k$ CISTs of $G_{n+1}.$
		\end{cor}
		
		According to the following definition, vertices situated at the center, referred to as central vertices, minimize the maximum distance to all other vertices in the graph.
		
		\begin{definition}[\cite{we}]
			If $G$ has a $u, v$-path, then the {\it distance} from $u$ to $v$, written $d_G(u, v)$ or simply $d(u, v),$ is the least length of a $u, v$-path. The {\it eccentricity} of a vertex $u,$ written $\epsilon{(u)},$ is $\displaystyle max_{v \in V(G)}d(u, v).$	The {\it center} of a graph $G$ is the set of all vertices of minimum eccentricity.  
		\end{definition}
		
	We will now present a construction of three completely independent spanning trees (CISTs) in the hypercube \( Q_n \) for \( n \geq 8 \), employing an inductive approach. Two of these trees will have diameters of at most \( 2n + 1 \) and \( 2n + 3 \), respectively, while the third tree will have a diameter of at most \( 2n + 4 \).
	
	\begin{thm}
		Let \( n \geq 7 \) be an integer. In the hypercube \( Q_n \), there exist three completely independent spanning trees such that two of these trees have diameters of at most \( 2n + 1 \) and \( 2n + 3 \), respectively, and the third tree has a diameter of at most \( 2n + 4 \).
	\end{thm}
		\begin{proof} 
			
			We already showed that there are three CISTs in \( Q_7 \) in Section 3. By using the CISTs \( T_1 \), \( T_2 \), and \( T_3 \) from \( Q_7 \) and applying Corollary 4.1, we can create three CISTs in \( Q_n \) for \( n \geq 8 \) using an inductive approach. 
			
			To minimize communication delays by keeping vertices close together, we focus on the diameters of these trees. Specifically, \( T_1 \) and \( T_2 \) in \( Q_7 \) have diameters of 15 and 18, respectively, while \( T_3 \) has a diameter of 17. We will first address the construction of three CISTs in \( Q_8 \). 
			
			In \( Q_8 \), which is composed of \( Q_7^0 \), \( Q_7^1 \), and \( E_8 \), let \( T_1^0, T_2^0,\) and \( T_3^0 \) be the CISTs in \( Q_7^0 \) and \( T_1^1, T_2^1, \) and \( T_3^1 \) be their counterparts in \( Q_7^1 \). We only need to demonstrate the construction of one such tree. 
			
			Consider the tree \( T_1 \) in \( Q_7 \) with the longest path of length 15 (from node $37$ to node $101$) and central vertex \( 0010101(21) \). By prefixing this vertex with 0 and 1, we obtain central vertices \( u_1 = 00010101 \) and \( v_1 = 10010101 \) for \( T_1^0 \) and \( T_1^1 \), respectively. Define \( \overline{T_1} = T_1^0 \cup T_1^1 \cup \langle u_1, v_1 \rangle \). This tree \( \overline{T_1} \) will be a spanning tree with a diameter of 17 in \( Q_8 \). Similarly, by constructing spanning trees \( \overline{T_2} \) and \( \overline{T_3} \) with diameters of at most 20 and 19, respectively, in \( Q_8 \), we obtain CISTs \( \overline{T_1} \), \( \overline{T_2} \), and \( \overline{T_3} \) in \( Q_8 \) with diameters of at most 17, 20, and 19, respectively, as per Corollary 4.1.
			
			Consider now \( Q_{n+1} \) for \( n \geq 7 \), which decomposes into two hypercubes, \( Q_n^0 \) and \( Q_n^1 \), with vertex sets \( \{x^0_i : 1 \leq i \leq 2^n\} \) and \( \{x^1_i : 1 \leq i \leq 2^n\} \), respectively. Let \( T^0_1, T^0_2,\) and \( T^0_3 \) be the CISTs in \( Q_n^0 \) with diameters at most \( 2n+1, 2n+4, \) and \( 2n+3 \), respectively. The corresponding CISTs in \( Q_n^1 \) will be denoted as \( T^1_1, T^1_2,\) and \( T^1_3 \).
			
			We will use induction to prove that the diameters of \( \overline{T}_1, \overline{T}_2,\) and \( \overline{T}_3 \) in \( Q_{n+1} \) are at most \( 2n+3, 2n+6, \) and \( 2n+5 \), respectively. It suffices to demonstrate this for one tree, say \( T^0_1 \).
			
			Let \( P^0_1 = x^0_1 - x^0_2 - \cdots - x^0_{n-1} - x^0_n - \cdots - x^0_{2n+2} \) be the longest path in \( T^0_1 \), with a length of \( 2n+1 \). Let \( P^1_1 = x^1_1 - x^1_2 - \cdots - x^1_{n-1} - x^1_n - \cdots - x^1_{2n+2} \) be the corresponding path in \( T^1_1 \).
			
			The vertex \( x^0_{n+2} \) is near the center of \( P^0_1 \); thus, any vertex in \( T^0_1 \) is within a distance of \( n+1 \) from \( x^0_{n+2} \). Similarly, any vertex in \( T^1_1 \) is within a distance of \( n+1 \) from \( x^1_{n+2} \). 
			
			Construct \( \overline{T}_1 \) by combining \( T^0_1 \) and \( T^1_1 \) and adding the edge \( \langle x^0_{n+2}, x^1_{n+2} \rangle \) in \( Q_{n+1} \). The resulting spanning tree \( \overline{T}_1 \) will have a diameter of \( 2n+3 \). Since \( \overline{T}_1 \) includes all vertices from \( T^0_1 \) and \( T^1_1 \), the maximum distance between any two vertices \( x^0_i \) and \( x^1_j \) in \( \overline{T}_1 \) will be \( d_{T^0_1}(x^0_i, x^0_{n+2}) + 1 + d_{T^1_1}(x^1_{n+2}, x^1_j) \leq (n+1) + 1 + (n+1) = 2n+3 \).
			
			By applying the same construction method, we obtain CISTs \( \overline{T}_2 \) and \( \overline{T}_3 \) in \( Q_{n+1} \) with diameters of \( 2n+6 \) and \( 2n+5 \), respectively. 
			
		\end{proof}
	
	\vspace*{-1 cm}

		\section{Conclusion}
		
		In this paper, we present the following results:
		
		\begin{enumerate}
			\item We provide a necessary condition for \(k\)-regular, \(k\)-connected bipartite graphs to have \(\left\lfloor \frac{k}{2} \right\rfloor\) completely independent spanning trees (CISTs).
			
			\item We prove that for any positive even integer \(n\) with \(2 < n \leq 10^7 \), except when \(n = 2^r\) and \( n \in \{161038, 215326, 2568226, 3020626, 7866046, 9115426 \} \) there do not exist \(\frac{n}{2}\) CISTs in \(Q_n\).
			
			\item Additionally, we provide a method to construct three completely independent spanning trees in the hypercube \(Q_n\) for \(n \geq 7\), showing that Hasunuma's conjecture does not hold for the even integer \(n = 6\), but it does hold for the odd integer \(n = 7\).
		\end{enumerate}

		While the necessary condition proved in this paper does not guarantee that every \( k \)-regular, \( k \)-connected bipartite graph will admit \( \left\lfloor \frac{k}{2} \right\rfloor \) CISTs.\\
		
		In future, we aim to address the following open questions:

		\begin{enumerate}
			\item  Since Hasunuma's conjecture does not hold for even integer \(n = 6\), but it holds for the odd integer \(n = 7\), we are motivated to explore whether this holds true for all odd integers.
			
			\item Does the $n$-dimensional hypercube \(Q_n\) for $n = 2^r, r\geq 3$ admit \(\frac{n}{2}\) CISTs? 
			
		\end{enumerate}

		\section*{Acknowledgment}
		
		The second author gratefully acknowledges the Department of Science and Technology, New Delhi, India, for awarding the Women Scientist Scheme (DST/WOS-A/PM-14/2021(G)) for research in Basic/Applied Sciences.

\end{document}